\documentclass[12pt]{amsart}
\usepackage[ansinew]{inputenc}
\usepackage{amsfonts,amsmath}
\input epsf
\usepackage{a4wide}
\def\C{\mathbb{C}}
\def\D{\mathbb{D}}
\def\R{\mathbb{R}}

 \def\UU{{\mathcal U}}
\def\II{{\mathcal I}}      \def\JJ{{\mathcal J}}    \def\KK{{\mathcal K}}
\def\RR{{\mathcal R}}
\def\DD{{\mathcal D}}     \def\NN{{\mathcal N}}
\def\cal{\mathcal }

\newtheorem{Th}{Theorem}[section]
\newcommand{\bt}{\begin{Th}}
\newcommand{\et}{\end{Th}}
\newcommand{\etd}{$\square$\end{Th}}

\newtheorem{theorem}[Th]{Theorem}

\newtheorem{proposition}[Th]{Proposition}

\newtheorem{lemma}[Th]{{\bf Lemma}}

\newtheorem{corollary}[Th]{{\bf Corollary}}

\newtheorem{Proposicao}[Th]{Proposition}
\newcommand{\bprop}{\begin{Proposicao}}
\newcommand{\eprop}{\end{Proposicao}}
\newcommand{\epropd}{$\square$\end{Proposicao}}

\newtheorem{Lema}[Th]{Lemma}
\newcommand{\bl}{\begin{Lema}}
\newcommand{\el}{\end{Lema}}
\newcommand{\eld}{$\square$\end{Lema}}

\newtheorem{Corol}[Th]{Corollary}
\newcommand{\bcor}{\begin{Corol}}
\newcommand{\ecor}{\end{Corol}}
\newcommand{\ecord}{$\square$\end{Lema}}

\newtheorem{Defini}[Th]{Definition}
\newcommand{\bd}{\begin{Defini}}
\newcommand{\ed}{\end{Defini}}

\newtheorem{Claim}[Th]{Claim}
\newcommand{\bclaim}{\begin{Claim}}
\newcommand{\eclaim}{\end{Claim}}

\newtheorem{remark}[Th]{Remark}

\newcommand{\cqd}{\nolinebreak
\hfill $\Box$}

\title{{\bf Topological obstructions to smoothness for infinitely renormalizable maps of the disk}}

\author{\it F. J. Moreira}
\date{\today}

\begin{document}

\epsfclipon

 {\abstract{We analyze the signature type of a cascade
of periodic orbits associated to period doubling renormalizable
maps of the two dimensional disk. The signature is a sequence of
rational numbers which describes how periodic orbits turn each
other and is invariant by topological conjugacies that preserve
orientation. We prove that in the class of area contracting maps
the signature cannot be a monotone sequence. This explains why
classical examples of infinitely renormalizable maps due to Bowen,
Franks and Young cannot be achieved by smooth dissipative maps
showing that there are topological obstructions to realize
infinitely renormalizable maps in the area contracting case. }}

\maketitle
\section{Introduction}

In 1975, Bowen and Franks \cite{[BF]} gave the first example of a
$C^1$ diffeomorphism of the sphere $S^2$ which is Kupka-Smale and
possesses neither sinks nor sources. Using some more techniques
Franks and Young \cite{[FY]} improved this result in 1980  and got
a $C^2$ example. We will refer to these two examples as {\bf BFY}
models. Finally in 1989, a $C^\infty$ example was found
 by Gambaudo, Strien and Tresser \cite{[GST]} which we will refer as
{\bf GST}. It is still an open question wether there exists a real
analytic example in $S^2$ with the mentioned properties.

The three cited examples are obtained by finding first a
Kupka-Smale embedding of the two-disk with neither sinks nor
sources, then  glueing this embedding with its inverse to obtain
the desired diffeomorphisms of $S^2$. It turns out that the
embeddings of the two-disk are infinitely renormalizable by a
period doubling cascade  of disks surrounding a cascade of
periodic points of saddle type.

Following \cite{[GSuT]},   we can associate to a cascade of
periodic orbits,  a signature consisting of a sequence
$(\ell_n)_{n \geq 0}$ of rational numbers such that each $\ell_n$
describes how the orbits of period $2^{n+1}$ are linked to the
orbits of period $2^{n}$ (see next section for the details) and,
is invariant by orientation preserving topological conjugacies.
From the work of Gambaudo, Sullivan and Tresser, the signature is
a convergent sequence for $C^1$ maps. By simple computation we
easily derive in section \ref{Signatureforarea-contracting...}
that the signature of the type {\bf BFY} is a decreasing sequence
converging to~$0$.

The main result of this work deals with the obstruction to the
realization of a monotone signature in the class of area
contracting maps of the two dimensional disk. We prove (Theorem
\ref{mainobstruct}
 and its Corollary \ref{corolthesequencecannot...}) under mild
assumptions on the geometry of the periodic cascade that, in the
class of area contracting embeddings of the two disk (which
contains {\bf GST} models) monotone signatures cannot occur and
thus, the decreasing sequence  obtained in {\bf BFY} model cannot
be achieved for smooth dissipative maps.

The work is organized as follows. In section 2 we introduce the
concept of signature and related properties. Then, in subection
2.1 we describe the properties of infinitely renormalizable maps
needed to obtain the model types {\bf BFY} and {\bf GST} whose
construction is sketched in the two subsequent subsections 2.2 and
2.3. Afterwards, section 2.4 is devoted to the bounded geometry
property used throughout  this work. The main result and its proof
is established in section 3, where  the contents of subsection 3.2
plays a crucial role, since therein we reduce our problem to the
analysis of the signature of a multimodal endomorphism of the
interval, and allows us to use one dimensional techniques stated
in subsection 3.1.

\section{Cascades of periodic orbits} \label{topobstructions}
Let ${g}$ be an orientation-preserving homeomorphism of the
2-disk $\D^2$. A {\it cascade of periodic orbits } of ${g}$
is a sequence of periodic orbits $\{O_n\}_{n \geq 0}$ with
periods  $\{q_n\}_{n \geq 0}$ such that, for each $n \geq
1$, we have:
\begin{enumerate}
\item $q_n =a_n.q_{n-1}$ with $q_0
=1$ and $a_n >1$,
\item  there exists a collection of disjoint simple
closed curves  $  C_n^0, \ldots  C_n^{q_{n-1}-1},$ bounding the disjoint disks
$  D_n^0, \ldots  D_n^{q_{n-1}-1},$ with the following properties:
   \begin{enumerate}
      \item  each $D_n^i$ contains one point of $O_{n-1}$ and $p_n$ points of $O_n$,
       \item ${g}(C_n^i)$ is isotopic to $C_n^{i+1 {\rm  mod  }  (q_{n-1})}$ in the punctured disk $\displaystyle{D^2 \setminus \bigcup _{i \leq n} O_i}$,

     \item  $\displaystyle {\bigcup _{0 \leq i\leq q_{n}-1} D_{n+1}^i \subset \bigcup _{0 \leq i\leq q_{n-1}-1} D_{n}^i}$.
      \item  The diameters of the disks $ D_{n}^i$ go
to zero with $n$.
   \end{enumerate}
\end{enumerate}
  Let $\{f_t\}_{t\in [0,1]}$ be an arc of embeddings
joining the identity map ${f}_0$ to ${g} = f_1$, and $\{f_t\}_{t\in
\R}$ be the extended arc of embeddings joining the identity to
all iterates of ${g}$ defined by ${g}_t = f^{[t]} o f_{\{t\}}$
(where $[t]$ and $\{t\}$ denote the integer and decimal part of $t$,
respectively). To each cascade of periodic orbits  $\{O_n\}_{n \geq
0}$ we associate a {\it {signature}}  $s(\{O_n\}_{n \geq 0}) =
\{\lambda_n\}_{n \geq 1}$, where, for all $n \geq 1$,  $\lambda _n$
is a rational number, $\lambda_n =l_n(f)/q_n$ and $l_n(f)$ is an integer
defined as follows:  \noindent In one of the $D_n^i$'s, pick the
point $x_{n-1}$ of $O_{n-1}$ and a point $x_n$ of $O_n$. Then, $l_n(f)$
is the algebraic number of loops that the vector
$$f_t(x_n)-f_t(x_{n-1})\over ||f_t(x_n)-f_t(x_{n-1})||$$
 performs on the
unit circle when $t$ goes from 0 to $q_n$ (here $||.||$ stands for
the $\R^2$ norm). Clearly, this number $l_n$ is independent of the
choice of the $D_n^i$ and of the choice of the point $x_n$ in
$D_n^i$.

\begin{remark} Let $\{f'_t\}_{t\in [0,1]}$ be another arc
of homeomorphisms joining the identity map to ${g}$ and
$\{f'_t\}_{t\in \R}$ be the extended arc joining the
identity map to all the iterates of ${g}$.  We denote
$l'_n$ and  $\lambda'_n$ the quantities defined previously
but computed for  $\{f'_t\}_{t\in [0,1]}$. Then there
exists an integer $k$ such that  $l'_n = l_n + kq_n$ and
thus $\lambda'_n =\lambda_n + k$ for all $n\geq 1$.
\end{remark}
 If the limit of the $\lambda_n$'s exists when $n$ goes
to infinity, we call it the {\it {asymptotic rotation number}} of
the cascade $\{O_n\}_{n \geq 0}$ and denote it by  $\omega
(\{O_n\}_{n \geq 0}).$

\noindent
Let $\{O_n\}_{n \geq 0}$ be a cascade of periodic orbits of a map $f$ and $h$ an embedding of the 2-disk. The map $g=h^{-1}\circ f\circ h$ possesses also a cascade of periodic orbits $\{O'_n\}_{n \geq 0}$ given by $O'_n=h(O_n) $ for ${n \geq 0}$. An isotopy $\{f_t\}_{t \in  [0,1]}$ between the identity and $f$ gives rise, in an natural way, to an isotopy $(g_t= h^{-1}\circ f_t \circ h)_{t \in [0,1]}$  between the identity and $g$. Computing $l_n(g)$ for this "conjugated" isotopy we have that
\begin{eqnarray}
l_n(g)&=&\left\{\begin{array}{rl}
        l_n(f) & \hbox{if $h$ preserves orientation}\\
        -l_n(f) & \hbox{if $h$ reverses orientation}
      \end{array}\right. \label{lnbyconjug} \ .
\end{eqnarray}

\noindent
Therefore, the signature of a cascade of periodic
orbits is a topological invariant, that is to say it is invariant
under conjugacy by an orientation-preserving homeomorphism of the
two disk.

\noindent
 It is easy to check that, given a sequence of rational
numbers $(\lambda_n)_{n\geq 0}$, one can construct an orientation-preserving homeomorphism of the 2-disk with a cascade of periodic
orbits $\{O_n\}_{n \geq 0} $ such that $s(\{O_n\}_{n \geq 0})
=\{\lambda_n\}_{n \geq 1}$. However there are topological
obstructions to realize such a cascade for an orientation-preserving
diffeomorphism of the 2-disk:
 \begin{theorem} \cite{[GSuT]} Any
cascade of periodic orbits of a $C^1$ orientation-preserving
diffeomorphism of the 2-disk possesses an asymptotic rotation
number.
\end{theorem}

\subsection{Infinitely
renormalizable maps}

Consider the unit disk $\D^m$ and a continuous map $g:{\D^m}\to
{\D^m}$. Assume that there exists in $\D^m$ a topological disk
$D$, that is to say, the image of $\D^m$ by a $C^0$-embedding
$\xi(g):{\D^m}\to {\D^m}$ ($\xi({\D^m})=D$),  such that there
exists $p \geq 2$ satisfying:
\begin{enumerate}
\item   $ D$, $g(D)$, \ldots  $g^{p-1}(D)$  have disjoint
interiors, \noindent
\item ${g}^p (D) \subset D$.
\end{enumerate}
We call such a map \em  renormalizable \em and the disk $D$
a \em domain of renormalization. \em  In this setting it is natural
to define the map
\begin{eqnarray*}
 \RR(g)&=&\xi^{-1}(g)\circ g^p \circ \xi(g)
\end{eqnarray*}
which turns out to be the first return map in the disk $D$
rescaled to the disk $\D^m$. We call the map $\RR(g)$ \em the
renormalized map. \em

\noindent The functional operator $\RR$, which associates a
renormalized map to a renormalizable one, is called \em the
renormalization operator. \em

\noindent If the renormalized map $\RR(g)$ is again renormalizable
 we
say that $g$ is \em twice renormalizable.\em It follows that we can
define in a natural way \em $m$-times renormalizable maps \em and
also \em infinitely renormalizable maps. \em

More precisely, a continuous map of the unit disk to itself is
infinitely renormalizable if there exists a sequence ${\D^m}
\supset D_0  \supset D_1 \supset \dots \supset D_n \dots$ of nested
topological disks, and  a sequence $(a_n)_{n \geq 0}$ of integers
greater than or equal to 2,  such that, for each $n \geq 0$:
\begin{enumerate}
\item  $ D_n$, $g(D_n)$, \ldots $g^{a_o\cdots a_n-1}(D_n)$
have disjoint interiors,
\item ${g}^{a_0.a_1\dots a_n} (D_n) \subset D_n$.
\end{enumerate}
When  more precision is required, we shall say that such a map is
\em $(a_n)_{n \geq 0}$-infinitely renormalizable.\em  We say that the
nested sequence of disks $(D_n)_{n \geq 0}$ is a \em cascade \em for
$g$ and the sets $ f^i (\DD_n) $, for $0 \leq i\leq a_0.a_1\dots
 a_n -1$, are called  \em the atoms of generation $n$ of ${g}$. \em
Denoting  by $\alpha_k(g):\D^m\to \D^m$ an
embedding whose
image is $D_k$ the \em $k$-renormalization of $g$ \em is defined as
\begin{eqnarray*}
\RR_k(g) & = & \alpha_k^{-1}(g) \circ g^{a_0 \cdots
a_{k-1}}\circ \alpha_k(g) \ \  .
\end{eqnarray*}

It is easy to check that the map $\RR_k(g)$ is $(b_n)_{n \geq
0}$-infinitely renormalizable with $b_n=a_{n+k}$ for each $n
\geq 0$ and the renormalized map of $\RR_{k}(g)$ is the map
$\RR_{k+1}(g)$, under the natural scaling
$\xi(\RR_{k}(g))=\alpha_{k}^{-1}\circ \alpha_{k+1}$, that is to say
$$ \RR_{k+1}(g) \ = \ \xi^{-1}(\RR_{k}(g))\circ \RR_k(g)^{a_k} \circ
\xi(\RR_{k}(g)) \ = \ \underbrace{\RR\circ \cdots\circ\RR}_k(g) \ .
$$

We say that an infinite renormalizable map is  \em  of bounded
combinatorial type \em  if the sequence $(a_n)_{n\geq 0}$ is
bounded.

In the case of the interval ${\D^1}=[-1,1]$, it is natural
to consider only the change of variables that are affine
maps, since the image of any embedding $\alpha:[-1,1]\to
[-1,1]$ is a nondegenerate interval, and thus we can carry
$[-1,1]$ onto $\alpha([-1,1])$ by an affine map.

 Let $\mathcal U$ be the set of
real-analytic maps $g:[-1,1]\to[-1,1]$, satisfying
\begin{enumerate}
 \item $g$  is strictly increasing in $[-1,0]$ and strictly
decreasing in $[0,1]$, that is to say, $g$ is unimodal with $0$ as
the critical point.
 \item $g(0)=1$ and $g^{\prime\prime}(0)\not=0$.
 \item By setting $a=-g(1)$, $b=g(a)$ we have
      \begin{enumerate}
                                \item $0<a<b<1$.
        \item $g(b)=g^2(a)<a$.
      \end{enumerate}
\end{enumerate}

\noindent Any map $g\in \mathcal U$  maps $[-a,a]$ onto $[b,1]$,
and $[b,1]$ onto $[-a,g(b)]\subset [-a,a]$. Moreover $g\circ g $
is again unimodal in $[-a,a]$. Then, in $\mathcal U$ we can define
the renormalization operator $\RR$ by
\begin{eqnarray*}
\RR(g)(x) &=& {1 \over g(1)}g^2(g(1)x) \ .
\end{eqnarray*}

\noindent In order to explain   quantitative  universal phenomena
appearing in bifurcation in one parameter families of maps in the
class $\mathcal U$, Coullet and
 Tresser (\cite{[CT]}, \cite{[TC]}) and   Feigenbaum \cite {[Fei]},
conjectured the following scenario  for the structure of the
renormalization operator defined above.

\begin{proposition} \label{ctf} \label{phirho}
The operator $\RR$ is a bounded $C^2$
operator having a fixed point $\phi \in \cal U$ with the following
properties:
\begin{enumerate}
 \item $\phi(x)=r(x^2)$, where $r$ is an analytic homeomorphism
defined in a neighbourhood of $[-1,1]$ in $\C$. In particular
$\phi$ is symmetric.
 \item The derivative $D\RR(\phi)$ is a compact operator whose
spectrum has a unique eigenvalue $\delta=4.6692...$ outside the
unit disk and all the other eigenvalues have modulus less than $1$. Let $\rho$ be the eigenvector corresponding to the eigenvalue $\delta$. This eigenvector is of the form \linebreak $\rho(x)=v(x^2)$,  where $v$ is an analytic homeomorphism defined in a neighbourhood of $[-1,1]$ in $\C$.
\end{enumerate}
\end{proposition}
Lanford in 1982 \cite{[La1]} gave the first complete proof of the
above proposition. This  proof is a rigorous analysis of the renormalization operator, based on numerical estimations made by computer. Later,
in 1987  Eckman and Wittwer \cite{[EW]} gave a different proof of
the same conjecture. Before that, the existence of fixed point for
the renormalization without the characterization given in \ref{ctf}
had already be proved by Campanino and Epstein \cite{[CaE]}.

\noindent
Using hyperbolic theory, the spectral property of the renormalization
operator given in Proposition \ref{ctf}, gives locally,
 the following picture:
\begin{proposition} In a neighbourhood of the fixed point $\Phi$ in
the space $\cal U$, there exists an unstable manifold $W^u$ of
dimension one tangent to an eigenvector $\rho:[-1,1]\to[-1,1]$
associated to de eigenvalue $\delta$ and a codimension one stable
manifold $W^s$ where for every map $g \in W^s$, $\RR_n(g)$ converges
exponential fast to $\Phi$ as $n \to\infty$. \end{proposition}

\noindent
In 1992 Sullivan \cite{[Su]} gave a conceptual explanation of
universality in renormalization, opening new directions in the Theory of Dynamical Systems.

\bt \cite{[Su]}
 There exists a space $\cal E$ of maps $g:[-1,1]\to [-1,1]$ of the
form
$${g}= \psi_1 \circ Q \circ\psi_2 \ ,$$
where $Q(x)=x^2$,  $\psi_1,\psi_2:[-1,1]\to [-1,1]$ are  orientation-preserving
and orientation-reversing  respectively, with the following properties
\begin{enumerate}
\item Every $f \in \cal E$  has an  holomorphic extension
$f^{\C}:W\to \C$ to a neighbourhood $W$
of $[-1,1]$ in $\C$.
\item $\cal E$ has a strongly compact subset  $\cal C$, that is to
say, for every sequence $(f_n)_{n \geq 1}$ in $\cal C$, the
holomorphic extension $(f_n^{\C})_{n \geq 1}$ has an uniformly
converging subsequence.
\item For any map $f \in \cal E$ which is infinitely renormalizable
and of the bounded combinatorial type there exists $n_0(f)$, such
that $\RR_n(f)\in \cal C$ for every $n \geq n_0$.
\item $f,g \in \cal C$ have the same combinatorial type if and only
if $(\RR_n^{\C}(f)-(\RR_n^{\C}(g))_{n \geq 0}$ converges uniformly to
$0$.
\item For any $C^{k+Lip}$ ($k \geq 1$)\footnote{Actually
$C^{k+Zygmund}$. See for example \cite{[MS]}} unimodal map
$f:[-1,1]\to [-1,1]$ $(a_n)_{n \geq 0}$-infinitely renormalizable
with bounded combinatorial type,  the $C^1$ distance of $\RR_n(f)$
to $\cal C$ goes to $0$ as $n$ goes to $\infty$. Moreover if
$(a_n)_{n \geq 0}$ is periodic of period $N$, then $(\RR_{nN}(f))_{n
\geq 0}(f)$ converges in the $C^k$-topology to a point
$\phi_{(a_n)_{n \geq 0}} \in \cal C $ that is a periodic point of
the renormalization operator, $\RR_N(\phi_{(a_n)_{n \geq
0}})=\phi_{(a_n)_{n \geq 0}}$. \end{enumerate}   \et

\noindent
The techniques introduced in the proof of this theorem are a
beautiful combination of real and complex analysis. Let us
finish this section pointing out a step in Sullivan's proof often
referred to  as \em "real bounds." \em

\begin{theorem}\label{su}  {\bf {\cite{[Su]}}} Let $f$ be an
$C^{k+Lip}$ ($ k \geq 1$)  infinitely renormalizable unimodal
map with  combinatorial type bounded by
 $N$. Then, for all $n\geq 0$:
\begin{enumerate}
\item The renormalized maps $\RR^n(f)$  have $C^i$-norm
bounded by a constant $L_f(i)$ which depends only on ${f}$ and $i$
($1\leq i\leq k$). (The Lipschitz constant of the $k$- derivative
is also bounded by a constant  depending only on $f$)
\item There
exist two constants $a_f$ and $b_f$  which depend only on
 ${f}$ such that, if ${\cal I}$ is an atom of the generation $m$  and ${\cal J} \subset {\cal I}$ is an atom of the
 generation $m+1$ of
 $\RR^n(f)$ then,
$$0< a_f \leq  {|{\cal J}| \over |{\cal I}| } \leq b_f <1$$
(where $|.|$ stands for the diameter).
\item All these bounds are ``beau" (bounded and eventually
universally
 (bounded)), that is to say, that for $n$ big enough, these bounds
 can be chosen so that they depend only on $N$.
\end{enumerate}
\end{theorem}

\noindent
In condition 2 the same estimate is true if ${\cal J}$ is a
connected component of the complementary of the atoms of generation
$m+1$ (a "gap") in ${\cal I}$. So for two atoms ${\cal J},{\cal K}
\subset {\cal I}$ of generation $m+1$ we  also have ${\rm dist}({\cal
J},{\cal K}) \geq a_f |{\cal I}|$.

\noindent
In dimension 2,
 infinitely renormalizable  maps are also frequently observed. For instance, they
appear naturally  in the  infinitely dissipative situation for a map
  $(x,y) \mapsto (g(x),0),$ where $g$ is an
 infinitely renormalisable map on the interval, and also in the  area preserving
case
 of a map exhibiting resonant islands.

\noindent
In the sequel we give two examples of  $(2)_{n\geq 0}$-infinitely
 renormalizable embedding of the two-disk.

\subsection{The BFY Model} \label{bowenfranks}

This construction was first introduced by R. Bowen and J. Franks
\cite{[BF]} and provided a first example of $C^1$ Kupka-Smale
diffeomorphism of the 2-sphere with no sinks nor sources. Later it
was improved to a $C^2$ example by Franks and Young \cite{[FY]}.

\noindent
Consider an  orientation-preserving  $C^\infty$ diffeomorphism  $g$ of
the unit disk, which satisfies the following properties:
\begin{enumerate}
\item ${g}$ is the identity in a thin annulus $A_\epsilon = \{x\in
\D, \Vert x\Vert \geq 1-\epsilon\}$;
\item there exist four disjoint disks $D(0)$, $D(1)$, $D(2)$,   and
$D(3)$ of the same radius $\rho_0 >{ 1 \over 4}$ such that
${g}$ restricted to an $\epsilon$-neighborhood $D_{\epsilon}(i)$ of
each disk $D(i)$ is a translation which maps the disk $D(i)$ exactly
onto the disk $D(i+1\ {\rm  mod}\ 4)$ and $g^4$ restricted to $D(0)$
is the identity.
\item In the interior of $\D\setminus (D(0)\cup D(1)\cup D(2) \cup
D(3))$, ${g}$ is a Morse-Smale system.
\end{enumerate}

\begin{figure}[h]
\epsfysize=6.2cm \centerline{\epsffile{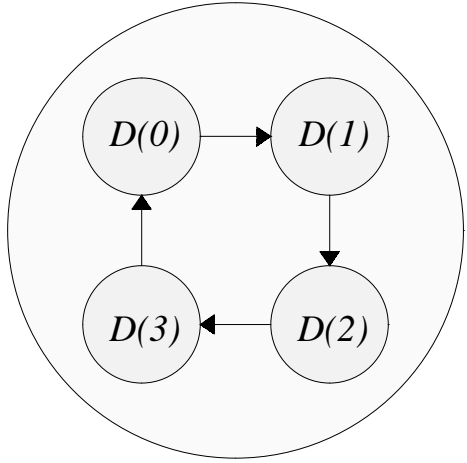}} \caption{
The map $g$} \label{gini}
\end{figure}

\noindent
Consider now an isotopy $G : [0,1]\times \D\to \D$ joining the
identity, $id=G(0,\cdot )$ to $G(1,\cdot )=g$ such that, for all
$t\in [0,1]$,   $G(t,\cdot )$ restricted to $A_\epsilon$ is the
identity and it maps rigidly the disks $D(i)$. For each positive
integer $N$ and $0 \leq i \leq N-1$, by setting
$$ g (N,i) \ =\
G({i+1 \over N}, \cdot)\circ G^{-1}({i \over N}, \cdot)
$$
we can write $g$  as the composition of $N$ maps $g(N,i)$,
$$
g \ = \ g(N,N-1)\circ \cdots \circ g(N,1) \circ g(N,0)
$$
that are $C^r$-close to the identity as $N$ goes to infinity,
$$
||g(N,i)-id||_r \leq \ { K(G,r) \over N} \ .
$$

\noindent
Let $\Lambda$ be an affine map that carries $\D$ to $D(0)$. Now we
are going to construct a sequence of functions $(f_n)_{n \geq 0}$,
($f_0=g$), $(4)_{1 \leq i \leq n+1}$-renormalizable with
$\RR_i(f_n)= f_{n-1}$ for $1\leq i\leq n$. Moreover, $f_n$  has
$4^{n+1}$ disks given by   $D_n(i)=f_n^i(\Lambda^{n}(\D))$, $(0 \leq
i \leq 4^{n+1}-1)$, where $f_n$ is a rigid move and $f_n^{4^{n+1}}$
restricted to $D_n(0)$ is the identity (this means that we can
define $\RR_{n+1}(f_n)$ as  the identity). The construction follows
by induction. Define $f_1$ as

$f_1(x)   \ =\ g(x)$ if $x$ is not in some $D(i)$, ($0 \leq i\leq 3$)

$f_1(x)   \ =\ g^{i+1}\circ\Lambda\circ g(4,i)\circ\Lambda^{-1}\circ\ g^{-i}(x)$ if $x$ is in some $D(i)$.

\noindent So, we have by the change of variables $\xi(f)=\Lambda$,
${\cal  R}(f_1)=  \Lambda^{-1}\circ f_1^4\circ\Lambda = g$. This
function has the $16$ disks with the properties needed for the first
step of induction (recall that $G(t,\cdot)$ and thus, $g(k,i)$ rigidly
permutes the four disks $D(0)$, $D(1)$, $D(2)$, $D(3)$).

\begin{figure}[h]
\label{f_2def}

\epsfysize=6.2cm
\centerline{\epsffile{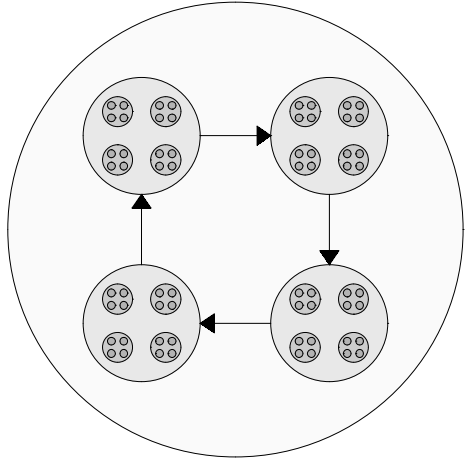}} \caption{
The atoms for the map $f_2$}
\end{figure}

\noindent
 Generally we define recursively $f_{n}$ by

$f_n(x)   \ =\ f_{n-1}(x)$ if $x$ is not in some $D_{n-1}(i)$, ($0 \leq i\leq 4^n-1$)

$f_n(x)   \ =\ f_{n-1}^{i+1}\circ\Lambda\circ g(4^n,i)\circ\Lambda^{-1}\circ\ f_{n-1}^{-i}(x)$ if $x$ is in some $D_{n-1}(i)$.

\noindent The map $f_n$ is again a $C^\infty$ diffeomorphism such that
$\RR(f_n)=f_{n-1}$ with $\xi(\RR(f_n))=\Lambda$. The annulus
$A_\epsilon$ and its reduced copies allows the surgery which
transforms $f_{n-1}$ into $f_n$ to be arbitrarily smooth.

\noindent
Since $f_n$ differs from $f_{n-1}$ only in the disks $D_{n-1}(i)$, where $f_{n-1}$ is a  rigid map, we have that
\begin{eqnarray*}
||f_{n}-f_{n-1}||_r &=& \max_{i}||\Lambda\circ g(4^n,i)\circ\Lambda^{-1}-id||_r \\
 & \leq & \max_i ||g(4^n,i)-id||\cdot \rho_0^{1-r}\\
 & \leq\ & K(G,r)\cdot {\rho_0^{(1-r)n}\over 4^n} \ .
\end{eqnarray*}
So, $f_n$ will be a converging sequence in the $C^r$ topology if we ensure that
\begin{eqnarray}
\sum_{n =1}^\infty{K(G,r)\cdot {\rho_0^{(1-r)n}\over 4^n}} &<&\infty
\label{serie} \end{eqnarray}
Thus, by the choice $\rho_0 > {1 \over 4}$ this is true for $r \leq
2$. The $C^2$ limit sequence $f_\infty$ is $(4)_{n \geq 1}$-
infinitely renormalizable by the scaling map $\Lambda$. Indeed this
map is a fixed point of the renormalization operator $f\mapsto
\Lambda^{-1}\circ f^4\circ \Lambda$ for the maps of $\D$ cyclically
permuting the four disks $D(0)$, $D(1)$, $D(2)$, $D(3)$.  This map
$f_\infty$ is $C^{2+\epsilon}$ but not $C^3$ (see \cite{[G]}). The
map $f^2_\infty$  obviously is a fixed point for the period doubling
operator $f\mapsto\Lambda f^2 \Lambda^{-1}$.

\noindent Starting with an map that rigidly permutes $a_0$ disks of radius $\rho_0> a_0^{-1}$ we can arrive, with this type of construction, to an $(a_0)_{n\geq 0}$-infinitely renormalizable map of the two disk. For any $a_0 >2$ we can find in $\D^2$  $a_0$ disks with radius $\rho_0 > a_0$. However, there is no space in $\D^2$ for two disks with radius $\rho_0>{1 \over 2}$. This is the reason why we construct an $(2)_{n \geq 0}$-infinitely renormalizable map in terms of an $(4)_{n \geq 0}$-infinitely renormalizable one.

\noindent The same kind of  techniques described in this section apply to higher dimensions $m \geq 3$, yielding, as explained by Gambaudo et Tresser in \cite{[GT2]}, $C^{m+\epsilon}$ infinitely renormalizable maps on $\D^m$.

\subsection{The GST model} \label{doublingmodel}

Now we do a brief presentation of an example of real analytic
embedding of the  2-disk which is infinitely renormalizable and
yields a $C^\infty$  Kupka-Smale diffeomorphism of the 2-sphere
with no sinks nor sources as described by J.M. Gambaudo, S. van
Strien and C. Tresser \cite{[GST]}. Unlike the Bowen-Franks model,
it does not result from a construction but from an exploration of
the properties of the period doubling  operator in one-dimension.

\noindent
Let $\phi(x)=r(x^2)$ and  $\rho(x)=v(x^2)$ be the maps given by Proposition  \ref{phirho}. Consider a small neighbourhood of $[-1,1]$ in $\R^2$
$$
R(\Delta)\ = \ \{(x,y)\in \R^2 \ | \ \exists x_0 \in [-1,1] \ {\rm with }  \ |x-x_0|+|y|<\Delta \}
$$
and, the corresponding neighbourhood in $\C^2$
$$
C(\Delta)\ = \ \{(x,y)\in \C^2 \ | \ \exists x_0 \in [-1,1] \ {\rm with }  \ |x-x_0|+|y|<\Delta \} \ .
$$
Let $H(\Delta)$ be the set of bounded analytic maps $G:C(\Delta):\to \C^2$ such that  $G(R(\Delta))\subset \R^2$.

\noindent Define the map $\Psi_\alpha$ by setting
\begin{eqnarray*}
\Psi_\alpha(x,y) &=& (r(x^2-\alpha y),0) \ .
\end{eqnarray*}
We know that $r$ is defined and is analytic in a neighbourhood of $[-1,1]$ in $\C$. Thus, for $\Delta$ small we have that the map  $(r(x),0)$ belongs to $H(\Delta)$. This implies, for $\Delta,\alpha>0$ sufficiently small that $\Psi_\alpha \in H(\Delta)$.

\noindent Since $\phi(x)=r(x^2)$ is a fixed point for the operator
$\RR(g)(x)= \lambda^{-1}g^2(\lambda x)$ ($\lambda = \phi(1)$) for
maps in $\cal U$, it is straightforward to verify that, by setting
$\Lambda(x,y)=(\lambda x,\lambda^2 y)$, the map $\Psi_\alpha$ is for $0<\alpha <1$
  a fixed point for the operator $\NN$ defined by:
\begin{eqnarray*}
\NN(g)&=&  \Lambda^{-1}\circ g^2\circ\Lambda \ .
\end{eqnarray*}

\noindent
Collet, Eckmann and Koch \cite{[CEK]} have  proved that $\NN$ is in fact a $C^2$ bounded operator from a neighbourhood of $\Psi_\alpha$ in $H(\Delta)$ into $H(\Delta)$. In the same work it is proved that $D\NN(\Psi_\alpha)$ is a compact operator from $H(\Delta)$ into itself. Although $1$ is an eigenvalue of $D\NN(\Psi_\alpha)$, in the cited work \cite{[CEK]} it is introduced another  operator $T_\alpha$ defined in a neighbourhood of $\Psi_\alpha$ in  $H(\Delta)$ into $H(\Delta)$, for which $\Psi_\alpha$ is an
hyperbolic fixed point. This operator $T_\alpha$ is a small perturbation of the operator $\NN$. We summarize the properties of $T_\alpha$ in the next proposition.

\begin{proposition}\label{psicek}  \cite{[CEK]} For  $\Delta,\alpha>0$ sufficiently small, there exists a $C^2$ transformation $T_\alpha$ from a neighbourhood of $\Psi_\alpha$ in $H(\Delta)$ into $H(\Delta)$ with the following properties.
\begin{enumerate}
\item $T_\alpha(\Psi_\alpha) \ = \ \Psi_\alpha$.
\item $T_\alpha(G)$ reads:
$$T_\alpha(G) \ = \ \Lambda_G^{-1}\circ G^2 \circ \Lambda_G \ ,$$
where $\Lambda_G:C(\Delta)\to C(\Delta)$ is an embedding belonging to $H(\Delta)$  such  that $\Lambda_G\to \Lambda$ as $G\to \Lambda$.
\item The derivative $DT_\alpha(\Psi_\alpha)$ of $T_\alpha$ in $\Psi_\alpha$ is a compact operator whose
spectrum has a unique eigenvalue $\delta=4.6692...$ outside the unit disk and all the other eigenvalues have modulus less than $1$. The eigenvector corresponding to this eigenvalue $\delta$ is of the form $(v(x^2-\alpha y),0)$.
\end{enumerate}
\end{proposition}

\noindent Let us set
$$B_G(\Delta) \ = \ \Lambda_G(R(\Delta)) \ .$$
Since $\Lambda_G$ sends $C(\Delta)$ into $C(\Delta)$ and $\Lambda_G \in H(\Delta)$ we have that $B_G(\Delta)\subset R(\Delta) $. Besides this,  from the fact that for  $\Delta >0$ small $\Phi([-1-\Delta,1+\Delta])\subset ]-1-\Delta,1+\Delta[ $,   we can easily derive (see \cite{[GST]} for further details)  that, for $\Delta, \alpha >0$ sufficiently small and  every $G \in H(\Delta)$ sufficiently close to $\Psi_\alpha$, one has that  $G(R(\Delta)) \subset R(\Delta)$.

\begin{figure}[!h]
\epsfxsize=11.5cm \centerline{\epsffile{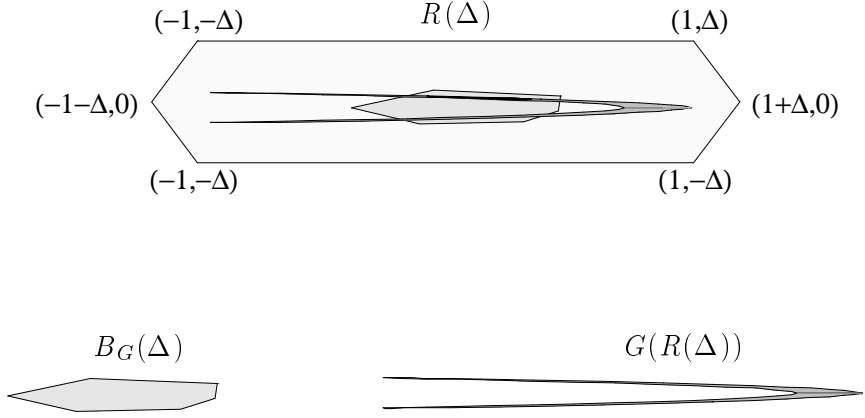}}
\caption{$G(R(\Delta))$ and $B_G(\Delta)$ for $G$ close to
$\Psi_\alpha$.}
\end{figure}

\noindent
From the  fact that $[\lambda,-\lambda] \cap \Phi([\lambda,-\lambda])=\emptyset$ and $\Phi^2[\lambda,-\lambda] \subset [\lambda,-\lambda]$,  it is simple to derive that for $\Delta,\alpha>0$ we have that $B_G(\Delta)$ and $G(B_G(\Delta))$ are disjoint and $G^2(B_G(\Delta)) \subset B_G(\Delta)$. This means that in fact $T_\alpha(G)$ is  a renormalized map of $G$.

\begin{figure}[!h]
\label{f(ax2-ay)} \epsfxsize=12.5cm
\centerline{\epsffile{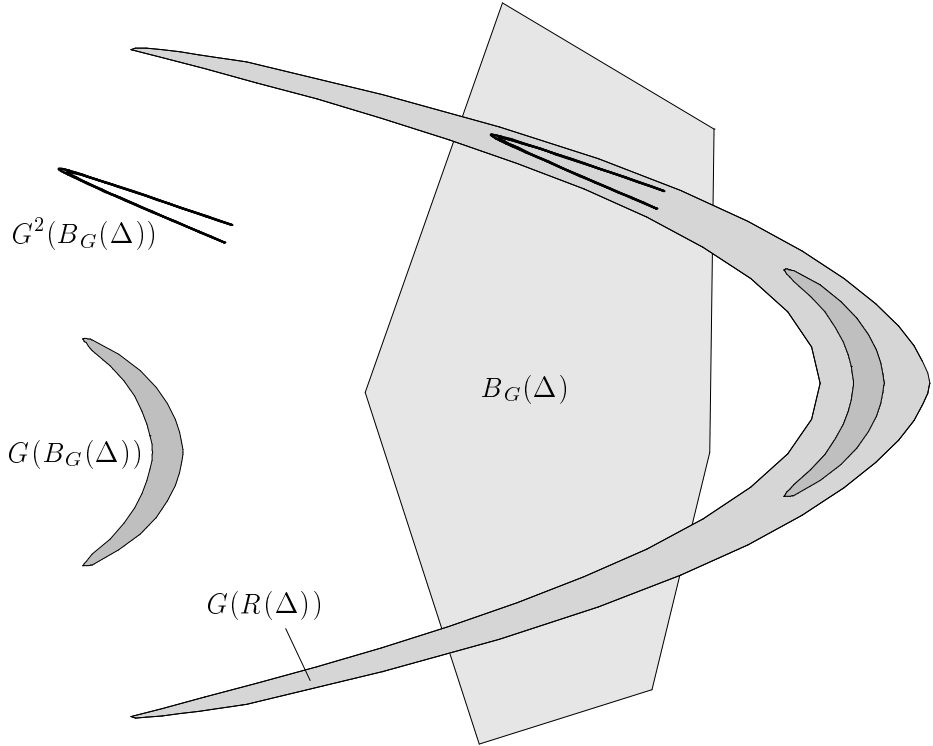}} \caption{$B_G(\Delta)$,
$G(B_G(\Delta))$ and $G^2(B_G(\Delta))$ for $G$ close to
$\Psi_\alpha$ (with a large zoom in the second coordinate).}
\end{figure}

\noindent Now, we know from Proposition \ref{psicek}, that
$\Psi_\alpha$ is a hyperbolic fixed point of $T_\alpha$.
Moreover, since $DT_\alpha$ has only one eigenvalue of modulus greater than one,
 the local stable manifold ${\cal W}_{loc}^s(T_\alpha)$ of $T_\alpha$
 in $H(\Delta)$ has codimension one. For every $G\in {\cal W}_{loc}^s(T_\alpha)$
 we have that $G$ is infinitely renormalizable\footnote{More precisely we should say the map $g_{|R(\Delta)}:R(\Delta)\to R(\Delta)$ is infinitely renormalizable} and its $n^{th}$ renormalized map is $T_\alpha^n(G)$.
 Consider a continuous family of real analytic embeddings $(f_{\mu,\epsilon})_{\mu \geq 0, \epsilon \geq 0}$ defined on $R(\Delta)$  satisfying the following conditions:

\begin{enumerate}
\item $f_{0,0} = \Psi_\alpha$
\item For every $\mu \geq 0$ and $\epsilon \geq 0$ $(f_{\mu,\epsilon})$  admits an extension in $H(\Delta)$.
\item $(\mu,\epsilon) \not= (\mu',\epsilon')$ implies that $f_{\mu,\epsilon}\not=f_{\mu',\epsilon'}$.
\end{enumerate}
For example, we can define $f_{\mu,\epsilon}(x,y)$ as
$$f_{\mu,\epsilon}(x,y)\ = \ (r(x^2-\alpha y)+\mu v(x^2-\alpha y),\epsilon x) \ .
$$
A family in these conditions  will necessarily cross , for small parameters $\mu$ and $\epsilon$, the  (codimension-one) local stable manifold ${\cal W}_{loc}^s(T_\alpha)$. Thus we have embeddings  $f_{\mu,\epsilon}$ infinitely renormalizable and their renormalizable maps $\RR_n(f_{\mu,\epsilon})=T^n_\alpha(f_{\mu,\epsilon})$ converge to the degenerate map $\Psi_\alpha$ as $n \to \infty$.

\subsection{Signature for maps with bounded geometry}

The maps found in \ref{bowenfranks} and \ref{doublingmodel} have a good behaviour under renormalization. They satisfy the following definition.

\bd   We say that an  infinitely renormalizable map $f$ of the $m$-disk, has $C^k$-bounded renormalization if it has bounded combinatorial type and  there exists a sequence of constants $(M_f(i))_{0 \leq i \leq k}$, where $M_f(i)$ depends only on ${f}$ and $i$ such that for all $n\geq 0$:
\begin{enumerate}
\item   the renormalized maps $\RR_n(f)$, the scaling maps  $\xi_n(f)=\xi(\RR_n(f))$,  and their inverse $\xi_n^{-1}$, are $C^{k+\epsilon}$ (in particular $f$ is $C^{k+\epsilon}$) and their $C^i$-norm  for $1 \leq i \leq k$  are bounded by $M_f(i)$.
\item $M_f(k)$ is also a  $\epsilon$-Holder constant  for $D^k\RR_n(f), D^k\xi_n(f)$ and $ D^k\xi_n^{-1}(f)$
\end{enumerate}
\ed

\noindent
Indeed in all the examples we have seen of infinitely renormalizable maps the scaling maps $\xi_n(f)$ are  $C^\infty$. They also are convergent as $n \to \infty$, and so they have the "beau" bound required in the above definition for any $k$. The differentiability of the bounded renormalization is given in all examples by the differentiability of the renormalized maps. In the case of a $C^{k+Lip}$ infinitely renormalizable unimodal  map of the interval, Theorem \ref{su} states that it has $C^k$-bounded renormalization. Moreover if such a map is $C^\infty$ it has $C^k$ bounded renormalization for every $k$. The Bowen Franks model in Section has $C^2$-bounded renormalization. Section \ref{doublingmodel} gives us an infinitely renormalizable map of the two disk with $C^k$-bounded renormalization for every $k$.

\noindent
Moreover all the atoms of this infinitely renormalizable maps have an universal bounded behaviour stated in the following definition.

\bd  We say that an infinitely renormalizable map $f$ of the $m$-disk, has $C^k$-bounded geometry if it has $C^k$-bounded renormalization and  there exists constants $0<a_f < b_f <1$  which depend only on
 ${f}$ such that for all $n\geq 0$, if ${\cal I}$ is an atom of the
 generation $m$ of $\RR_n(f)$ and ${\cal J} \subset {\cal I}$ is
 an atom of the generation $m+1$, then
  $$ a_f \leq\ {|{\cal J}| \over |{\cal I}|} \   \leq b_f \ ,$$
( where $|.|$ stands for the diameter). For another atom ${\cal K} \subset {\cal I}$ of generation $m+1$ we  also have
$$a_f \leq \ {{\rm dist}({\cal J},{\cal K})\over |{\cal I}|} \  \leq b_f  \ .$$
\ed

\begin{remark} Actually, the above relation in the ratios of the atoms forces ${f}$ to have bounded combinatorial type.
\end{remark}
\begin{remark} Bounded geometry implies that the diameter of an atom $I_l$ of generation $l$ goes to zero exponentially fast when  $l$ increases. In fact, it is straightforward to derive inductively that,
$$a_f^{l-1}\cdot|\D^m| \ <  \ |I_l| < \ b_f^{l-1}\cdot|\D^m| \ .$$
\end{remark}

\noindent
For simplicity we always say that a $C^k$ infinitely renormalizable map satisfies the bounded geometry hypothesis (respectively bounded renormalization hypothesis) if it satisfies the  $C^k$-bounded geometry properties (respectively $C^k$- bounded geometry  renormalization). If $f$ is $C^\infty$ we say that it has bounded geometry (respectively bounded renormalization) if it satisfies the $C^k$-bounded geometry properties (respectively the $C^k$-bounded renormalization properties ) for every $ k\geq 0$.

\noindent
 {\bf Remark:} To have bounded geometry is a strong
 assumption. An
 infinitely renormalizable map of the interval ${f} \in {\cal U}^{1+Lip}$ with
  bounded combinatorial type, satisfies this assumption (Theorem \ref{su}),
 and recently it has been proved that this is also the case
for other one-dimensional
 maps with finitely many critical points (see \cite{[Sm]}). However there is
 no result of this type for two-dimensional maps.

If  an $(a_n)_n$-infinitely renormalizable map $f$ has the
property that  any two atoms of the same generation are
disjoint (we recall that in our definition they  only have
disjoint interiors), it follows, from Brouwer fixed point
theorem, that the map ${f}$ possesses a sequence of
periodic orbits $\{O_n\}_{n \geq 0}$ with periods
$\{q_n\}_{n \geq 0}$ such that $q_n = a_n.q_{n-1}$ with
$q_0 =1$.

\noindent From now on we will consider an infinitely renormalizable map of the two-disk with the bounded geometry hypothesis. In particular, we know that under this hypothesis we have that  the distance of two distinct atoms of the same generation is greater than $0$. Thus they are disjoint.
\noindent
The collection of disjoint simple closed
curves  $  C_n^0, \ldots  C_n^{q_{n-1}-1},$ bounding the disjoint
disks $  D_n^0, \ldots  D_n^{q_{n-1}-1},$ where $D_n^i = f^i (\DD_n(f))$ satisfies:
\begin{enumerate}
\item  each $D_n^i$ contains one point of $O_{n-1}$ and $a_n$ points of
$O_n$,
\item ${f}(C_n^i)=C_n^{i+1 {\rm  mod  }  (q_{n-1})}$
\item  $\displaystyle {\bigcup _{0 \leq i\leq q_{n}-1} D_{n+1}^i \subset
\bigcup _{0 \leq i\leq q_{n-1}-1} D_{n}^i}$.
\item the diameter of the disks $D_n$ is less than $b_f^{n-1}|\D^2|$. Thus, it goes to zero as $n$ goes to infinity.
\end{enumerate}

\noindent
Therefore, we have associated to an infinitely renormalizable map,  a cascade of periodic orbits.
\begin{remark} The cascade of periodic orbits associated to an infinitely renormalizable map is not necessarily unique.
\end{remark}

\noindent
This cascade for $f$ gives rise to a cascade of periodic orbits for $\RR(f)$, where $O_{n}(\RR(f))= \xi(f)^{-1}\circ f^{a_0} \circ \xi(f)(O_{n+1})$.

\noindent
Consider now the signature of $f$, $({l_n(f)\over a_0 \cdots a_n})_{n \geq 1}$ computed by  an isotopy $\{f_t\}_{t \in [0,1]}$ of the disk joining the identity $f_0=id$ to $f_1=f$.  The image  $f_t({\cal D}(f))$ of the domain of renormalization ${\cal D}(f))=\xi(f)(\D^2)$ is a disk, with $f_0({\cal D}(f))={\cal D}(f)$ and $f_{a_0}({\cal D}(f))\subset {\cal D}(f)$. We can find a continuous arc $(A_t:{\cal D}(f)\to \D^2)_{t \in [0,a_0]}$ of affine maps, with $A_0=A_{a_0}=id$ and $A_t\circ f_t({\cal D}(f))\subset {\cal D_0}(f)$. With this arc of affine maps, we can construct an isotopy $F_t=A_t\circ f_t$, $t \in [0,a_0]$  of ${\cal D}(f)$. For the computation of $l_{n+1}(f)$, we look for the number of loops that the vector
$$f_t(x_{n+1})-f_t(x_{n})\over ||f_t(x_{n+1})-f_t(x_{n})||$$
performs when $t$ goes from $0$ to $a_0 a_1\cdots a_n$. Denote by $l'_{n+1}(f)$ the number of loops that the vector
$$F_t(x_{n+1})-F_t(x_{n})\over ||F_t(x_{n+1})-F_t(x_{n})|| \ .$$
performs when $t$ goes form $0$ to $a_0 a_1 \cdots a_n$.  Since $A_t$ are affine maps,  we have  that $l'_{n+1}(f)=l_{n+1}(f)$.
Computing $l_n(\RR(f))$ for the isotopy $\xi^{-1}(f)\circ F_t\circ \xi(f)$ we have, attending to \ref{lnbyconjug},  that $l_n(\RR(f))=l'_{n+1}(f)$, if $\xi(f)$ preserves orientation and,  $l_n(\RR(f))=-l'_{n+1}(f)$ if $\xi(f)$ reverses orientation. Thus we have
\begin{eqnarray}
l_n(\RR(f))&=&\left\{\begin{array}{rl}
        l_{n+1}(f) & \hbox{if $\xi(f)$ preserves orientation}\\
        -l_{n+1}(f) & \hbox{if $\xi(f)$ reverses orientation.}
      \end{array}\right. \label{lnbyren}
\end{eqnarray}

\noindent
Of course in the above equalities we assume that  the computation for $
l_n(\RR(f))$ is made using the isotopy $\xi^{-1}(f)\circ A_t\circ f_t \circ \xi(f)$.
By setting $\alpha_k= \xi_0(f)\circ \cdots \circ \xi_{k-1}(f)$, the  renormalized map $\RR_k(f)$ has a  cascade $(O_n(\RR_k(f))_{n \geq 1}$  given by $O_n(\RR_k(f))=\alpha_k(O_{n+k}(f))$.  Inductively we have for these two cascades:
\begin{eqnarray}
l_n(\RR_k(f))&=&\left\{\begin{array}{rl}
        l_{n+k}(f) & \hbox{if $\alpha_k$ preserves orientation}\\
        -l_{n+k}(f) & \hbox{if $\alpha_k$ reverses orientation.}
      \end{array}\right. \label{lnRk}
\end{eqnarray}

The equality given in \ref{lnbyren} allows us to compute immediately the signature for a fixed point of the renormalization operator under a fixed change of coordinates $\Lambda$.   If for an integer $a \geq 2$ and an embedding $\Lambda$ of the two disk  we have $\RR(f)=\Lambda^{-1}\circ f^a \circ \Lambda =f$ then,   for $ n \geq 1$ we have that:
\begin{eqnarray}
l_n(f)&=&\left\{\begin{array}{rl}
        l_{n+1}(f) & \hbox{if $\Lambda$ preserves orientation}\\
        -l_{n+1}(f) & \hbox{if $\Lambda$ reverses orientation}
      \end{array}\right.  \ .
\end{eqnarray}
and thus,
    \begin{eqnarray*}
 l_n(f)&=&\left\{\begin{array}{rl}
        l_{1}(f) & \hbox{if $\Lambda$ preserves orientation}\\
        (-1)^{n-1}l_{1}(f) & \hbox{if $\Lambda$ reverses orientation}
      \end{array}\right.  \ .
\end{eqnarray*}

\noindent
If we always choose an isotopy  for which  $l_1=1$, we have for the Bowen-Franks model ($\Lambda$ preserves orientation):
\begin{eqnarray*}
s(\{O_n\}_{n \geq 0}) &=& \left({1  \over 2^n}\right)_{n
\geq 0}
\end{eqnarray*} and for the period doubling model
\begin{eqnarray*}
s(\{O_n\}_{n \geq 0}) &=& \left({(-1)^{n-1}  \over
2^n}\right)_{n \geq 0} \ .
 \end{eqnarray*}

\section{Main result}
\label{Signatureforarea-contracting...}

The main result of this chapter says that for smooth infinitely
renormalizable embeddings of the two disk, area contracting  and
satisfying the bounded geometry hypothesis, the signature cannot
be a monotonic sequence.

\begin{theorem} \label{mainobstruct}
 Let $f$ be an $(a_n)_{n}$-infinitely
renormalizable embedding of the 2-disk with $a_n=2$ for all $n$,
and $({\ell_n \over 2^n})_n $the signature of the cascade of
periodic orbits associated to $f$. Assume that $f$ is $C^\infty$, has bounded
geometry and contracts the area. Then the sequence
$({\ell_n \over 2^n})_n $  alternates, that is to
say, for each $N>0$ we can always find $n_1$ and $n_0$ great or
equal to $N$, such that
 \begin{eqnarray*} {\ell_{n_0} \over
2^{n_0}} < {\ell_{n_0+1} \over 2^{n_0+1}} & {\rm and} & {\ell_{n_1}
\over 2^{n_1}} > {\ell_{n_1+1} \over 2^{n_1+1}} \ .
 \end{eqnarray*}
\end{theorem}
 As a consequence we get

\begin{corollary} \label{corolthesequencecannot...}
The sequence $({1 \over 2^n})_n$ cannot be the signature of an
$(2)_n$-infinitely renormalizable area-contracting embedding with
bounded geometry.
\end{corollary}

\noindent Roughly speaking, the
Bowen-Franks model cannot be realized in the class of
area-contracting embeddings with bounded geometry.  The proof of Theorem \ref{mainobstruct} is given in the last section, but before we need to recall a nice property of multimodal maps.

All the remaining work is devoted to the proof of the Theorem
above. First, in section \ref{..formultimodalmaps...} we recall
the results of multimodal maps that we need to achieve the proof.
In section \ref{Rigorousdynamicsreduction} we derive the reduction
of the dynamics of the initial map to the dynamics of a multimodal
endomorphism of the interval that allows to conclude the proof in
the last section.

\subsection{Cascade for multimodal maps of the interval}
\label{..formultimodalmaps...}

Let $g$ be a  $(2)_{n \geq 0}$infinitely renormalizable multimodal map of the interval $I=[c_0,c_{q+1}]$,  with critical points $c_1 < \ldots <c_q$. Denote by $I_k$ the interval of monotonicity $[c_{k-1},c_k[$. As usual $(\DD_n(g))_n$ denotes a cascade for $g$. We assume that $g$ satisfies the bounded geometry assumption, so that  atoms of the same generation are disjoint and the diameter of the atoms of generation $n$ goes to zero as $n \to \infty$.

\noindent
Under these conditions, in each atom $D_n^i=g^i(\DD_n(g))$, we have a periodic point $x_{n}$ of period $2^n$ and two periodic orbits $x_{1,n+1}, x_{2,n+1}$ of period $2^{n+1}$, with $x_{1,n+1}<x_n<x_{2,n+1}$. It is usual to call $x_n$ a \em father \em  of $x_{i,n+1}$ ($i =1,2$). A sequence of periodic orbits $(O_n)_{n \geq 0}$,  where for every $n$   the period of $O_n$ is $2^n$ and every point in $O_n$ has two sons in $O_{n+1}$, is called  a \em period-doubling cascade  of orbits \em for $f$.

\noindent
For a sufficiently large generation the father and its sons belong to the same interval of monotonicity or are in two adjacent intervals $I_k, I_{k+1}$, and the father is always between two sons.

\begin{figure}[h]
\centerline{\epsffile{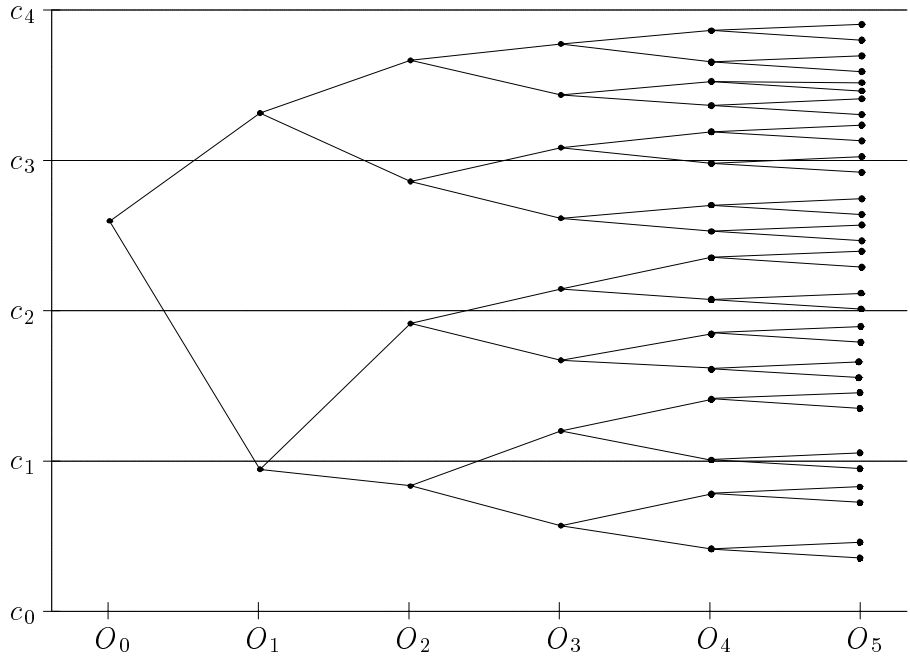}} \caption{The structure of
periodic orbits for multimodal maps. In this example, we have
$\Phi(1,3)=1$, $\Phi(2,3)=2$, $\Phi(3,3)=2$ and $\Phi(4,3)=3$}
\label{multimodalstructure}
\end{figure}

\noindent
Using the fact  that for any cascade of period doubling orbits  each father is between two sons  and, that globally we have the structure suggested by the Figure \ref{multimodalstructure}, Courcelle \cite{Olivier} was able to establish a useful relation for the number of periodic orbits in each interval of monotonicity:

\begin{proposition}  \label{oliv} Let $\Phi(k,n)$ denote the number of periodic orbits of period $2^n$ of a period doubling cascade  in the interval of monotonicity $I_k$ for $g$.
For all $n \geq 0$ and for all $k\in \{0, \ldots , q\}$
\begin{eqnarray*}  \Phi(k,n+1)&=&2\Phi(k,n)+r(k,n) \end{eqnarray*}
where for all $p_2 \geq p_1\geq 0$,
\begin{eqnarray*} \left|\sum_{j=p_1}^{p_2}{r(k,j)}\right| &\leq &2 \ . \end{eqnarray*}
\end{proposition}

\subsection{Rigorous dynamics reduction} \label{Rigorousdynamicsreduction}

All along this section we assume that $f$ is an
$(a_n)_{n\geq 0}$-infinitely renormalizable map   of the
two disk of the class $C^k$, where $k$ is  greater or equal
than $1$. \noindent We fix  a cascade $(D_n)_{n \geq 1}$ of
$f$ and its associated  sequence of embeddings
$(\alpha_n)_{n \geq 0}$ as defined in Section 1.1. So there
is no ambiguity on the definition of the scaling functions
$\xi_n(f)=\xi(\RR_n(f))=\alpha^{-1}_{n-1}(f)\circ\alpha_n(f)$.
The $n$-renormalized of $f$ is \begin{eqnarray*} \RR_n(f)
&=&  \alpha_n^{-1}(f)\circ f^{a_0\cdots a_{n-1}}\circ
\alpha_n(f) \\
         &=& \xi_{n-1}^{-1}(f)\circ\cdots\xi_0^{-1}(f) \circ f^{a_0\cdots a_{n-1}}
              \circ\xi_0(f)\cdots \xi_{n-1}(f) .
\end{eqnarray*}
Once  the scaling functions fixed, there is no ambiguity on the definition of the renormalization
operator  $\RR$ acting in the sequence of
renormalized maps $(\RR_n(f))_{n \geq 0}$:
$$
\RR(\RR_n(f))\ =\ \xi_n^{-1}\circ\RR_n^{a_n}\circ\xi_n \ = \ \RR_{n+1}(f) \ .
$$

\noindent We will also assume that $f$ satisfies the bounded
renormalization properties. We recall that under this assumption
$f$ is of the bounded  combinatorial type, that is to say $\{a_n,\
{n\geq 0}\}$ is a finite set, and thus, for  every subsequence
$(a_{\phi(n)})_{n \geq 1}$, there is some $p_0 \in \{a_n,\ {n\geq
0}\}$ with $a_{\phi(n)}=p_0$ for every $n \geq 1$. Moreover, the
bounded renormalization hypothesis  states that  the sequence of
renormalized maps $(\RR_n(f))_{n \geq 0}$ belongs to a space
${\cal U}^{k}={\cal U}^{k}(f)$  defined by the $C^{k+\epsilon}$
maps  $g:\D^2 \mapsto \D^2$, satisfying:
\begin{enumerate}
\item  $||D^{(i)}g|| \leq M_f(i)$ for every $ 1 \leq i \leq k$
\item The Holder constant $Hold(D^{(k)}g) \leq M_f(k)$.
\end{enumerate}

\noindent Thus, for every sequence $(g_n)_n$ in ${\cal U}^k$, we have that  for $i =0, \ldots , k$,
$(D^{(i)}g_n)_n$, satisfies the equicontinuity property:
 $\forall \epsilon >0 \ \exists \delta >0$ such that  for all $n \geq 0$
\begin{eqnarray*}
||x-y||<\delta \Rightarrow
||D^{(i)}g_n(x)-D^{(i)}g_n(y)|| < \epsilon \ .
\end{eqnarray*}

\noindent From the Arzela-Ascoli Theorem  there exists a
$C^k$-converging  subsequence of $(g_n)_n$. That is to say:  there
exists an increasing sequence of integers $(\phi(n))_{n \geq 1}$
and a $C^k$ map $g_\infty $  such that for every $0 \leq i \leq k$
the limit (in the $C^k$-topology) of the sequence
$(D^{(i)}(g_{\phi(n)}))_n$ exists and
$$
D^{(i)}(g_\infty) \ = \ \lim_{n \to \infty}D^{(i)}(g_{\phi(n)}) \ .
$$
By continuity, this map $g_\infty$ is also in ${\cal U}^k$.
To summarize:
\bt {\rm \bf (Arzela-Ascoli)} \label{ascoli-arzela}
For any sequence of maps $(g_n)_{n \geq 1}$ in $ {\cal U}^k$,  there exists a $C^k$-converging  subsequence  in $ {\cal U}^k$.
\et
\begin{remark}
This result implies that the limit of  any $C^0$-converging sequence in $\UU^k$ belongs to $\UU^k$. In particular such $C^0$-limit is of class $C^k$.
\end{remark}
When we refer to  an   accumulation point of a sequences in $\UU^k$ we mean that it is an accumulation point in the $C^0$ topology. By the previous remark, these $C^0$-accumulation points are also in $\UU^k$.
\begin{corollary} \label{notempty}
The sequence of renormalized maps $({\cal R}_n(f))_n$ has an accumulation point
${g_\infty}\in {\cal U}^{k}$.
\end{corollary}

\noindent Now let   $g$ be an arbitrary accumulation point of the sequence  $(\RR_n(f))_{n \geq 1}$ and ${\phi(n)}$ be a sequence such that $(\RR_{\phi(n)}(f))$ converges to ${g}$. From this sequence  we can extract  a subsequence ${\psi(n)}$ such that the sequence $(a_{\psi(n)})_n$ is constant (say equal to $p_0$) and the sequence ${\cal R}_{\psi(n)+1}(f)={\cal R}({\cal R}_{\psi(n)}(f))$  also converges in ${\cal U}^{k}$. Such a subsequence exists since $(a_n)_n$ takes only a finite number of values and ${\cal R}_{\phi(n)+1}(f)$ is a sequence in ${\cal U}^{k}$ (Arzela-Ascoli Theorem \ref{ascoli-arzela}). Besides this, we also can impose
 (once more, by the Arzela-Ascoli theorem) that the corresponding scaling maps $\xi_{\psi(n)}(f)$  also converge. Since $(\psi(n))_n$ is a subsequence of $(\phi(n))_n$, we still have $g =  \lim_{n\to \infty} {\cal R}_{\psi(n)}$.

\noindent Since for each $n$,
\begin{eqnarray*} \xi_{\psi(n)}(f)\circ {\cal R}_{\psi(n)+1}(f) &=& {\cal R}_{\psi(n)}^{p_0}(f)\circ \xi_{\psi(n)}(f) \ ,
\end{eqnarray*}
we obtain by continuity, taking  limits, the relation
\begin{eqnarray*}
\zeta_0\circ g &=& g^{p_0}\circ \zeta_0 \ ,
\end{eqnarray*}
where
\begin{eqnarray*}
g &=& \lim_{n\to \infty}{\cal R}_{\psi(n)+1}(f)\\
\zeta_0 &=& \lim_{n\to \infty}\xi_{\psi(n)}(f),
\end{eqnarray*}

\noindent From the bounded renormalization hypothesis, we also have  that for \linebreak $0 \leq i \leq k$,  $||D\xi^{-1}_n|| \leq M_f(i) $ \footnote{Notice that we cannot apply directly the Arzela-Ascoli Theorem to the sequence  $\xi_n^{-1}$ since its domain of definition depends on $n$.}. This says that for every $n$, $\xi^{-1}_n$  is at least Lipschitz  with Lipschitz constant smaller than $M_f(1)$, which means that:
\begin{eqnarray*}
||\xi^{-1}_n(f)(y_1) -\xi^{-1}_n(f)(y_2)|| & \leq & M_f(1)\cdot ||y_1-y_2||\ \ \forall y_1,y_2 \in \xi_n(\D^2) \ .
\end{eqnarray*}
This implies that for every $n$,
\begin{eqnarray*}
||\xi_{\psi(n)}(f)(x_1) -\xi_{\psi(n)}(f)(x_2)|| & \geq & { 1 \over M_f(1)}\cdot ||x_1-x_2||  \ \forall x_1,x_2 \in \D^2 \ ,
\end{eqnarray*}
which leads, after taking limits as $n \to \infty$,
\begin{eqnarray*}
||\zeta_0(x_1) -\zeta_0(x_2)|| & \geq & { 1 \over M_f(1)}\cdot ||x_1-x_2||  \ \forall x_1,x_2 \in \D^2 \ .
\end{eqnarray*}
From this last inequality it follows  that $\zeta_0$ is a homeomorphism over its image and its inverse $\zeta_0^{-1}$ has the same differentiability as $\zeta_0$.
From the fact that for each $n\geq 0$, we have:
$$
{\cal R}_{\psi(n)}^{p_0}(f)(\xi_{\psi(n)}(f)(\D^2))
 \subset \xi_{\psi(n)}(f)(\D^2).
$$
 and,  for $0\leq i <j\leq p_0-1$
\begin{eqnarray*}
{\rm Int}\left({\cal
R}_{\psi(n)}^{i}(f)(\xi_{\psi(n)}(f)(\D^2))\right)\bigcap
{\rm Int} \left({\cal
R}_{\psi(n)}^{j}(f)(\xi_{\psi(n)}(f)(\D^2))\right) & = &
\emptyset \ ,
\end{eqnarray*}
 we get,  by continuity:
$$ g^{p_0}({\zeta}_0(\D^2)) \subset {\zeta}_0(\D^2) \ , $$
and $$ {\rm Int} \left(g^i({\zeta}_0(\D^2))\right) \bigcap
{\rm Int} \left(g^j({\zeta_0}(\D^2)\right) \ = \emptyset\ .
$$ The above  induced statements  are summarized in the
following lemma. \bl \label{gren313} For any accumulation
point $g=\lim_{n\to \infty}{\cal R}_{\phi(n)}(f)\in {\cal
U}^k $ there exists  an accumulation point $p_0$  of the
sequence $(a_{\phi(n)})_n$, such that $g$  is
renormalizable by ${\cal R}(g) = \zeta_0^{-1}\circ
g^{p_0}\circ \zeta_0$. Furthermore, ${\cal R}(g)$ is a
limit point of the sequence $({\cal R}_{\phi(n)+1}(f))_n$
and $\zeta_0$ is an embedding which is an accumulation
point of the sequence $(\xi_{\phi(n)}(f))_n$. \cqd \el

\noindent Now  we are going to apply this lemma inductively. From an arbitrary accumulation point  $g$ of the sequence $({\cal R}_n(f))_n$,  we obtain a renormalization ${\cal R}(g)= \zeta_0^{-1}\circ g^{p_0}\circ \zeta_0$ which is  an accumulation point of the sequence  $({\cal R}_{\phi(n)+1}(f))_{n}$. Suppose that we have defined ${\cal R}_1 (g), \ldots , {\cal R}_m (g)$ with for $i =1, \ldots m$:
\begin{eqnarray*}
{\cal R}_{i}(g) &=& \zeta_{i-1}^{-1}\circ ({\cal R}_{i-1}(g))^{p_{i-1}} \circ \zeta_{i-1} \\
 & =& \left(\prod_{t=0}^{i-1} \zeta_t^{-1}\right)\circ g^{p_0\cdots p_{i-1}}\circ \left(\prod_{t=0}^{i-1} \zeta_t\right)
\end{eqnarray*}
where\footnote{Of course,   we are doing an incorrect use of the symbol $\prod$ in order to simplify some notation.}
\begin{enumerate}
\item $\zeta_i$ is an accumulation point of the sequence $(\xi_{\phi(n)+i-1}(f))_n$.
\item ${\cal R}_{i}(g)$ is a limit point of the sequence $({\cal R}_{\phi(n)+i}(f))_n$
\item $p_{i-1}$ is an accumulation point of the sequence $(a_{\phi(n)+i-1})\ . $
\end{enumerate}

\noindent Then we can apply the previous lemma to the accumulation point ${\cal R}_m (g) \in {\cal U}^{k}$ and obtain that ${\cal R}_m (g)$ is renormalizable by
$$
{\cal R}_{m+1}(g) \ = \ {\cal R}({\cal R}_{m}(g))  \ = \ \zeta_m^{-1}\circ ({\cal R}_m(g))^{p_m} \circ \zeta_{m}
$$
where
\begin{enumerate}
\item $\zeta_m$ is an accumulation point  of the sequence $(\xi_{\phi(n)+m})_n$,
\item ${\cal R}_{m+1}(g)$ is a limit point of the sequence  $({\cal R}_{\phi(n)+m+1})_n$,
\item $p_m$ is an accumulation point of the sequence $(a_{\phi(n)+m})_n\ .  $
\end{enumerate}
This yields:
\bl \label{Rm(g)}
Any accumulation point $g= \lim_{n\to \infty}{\cal R}_{\phi(n)}(f)$ is $(p_m)_{m \geq 0}$-infinitely  renormalizable with
\begin{eqnarray*}
{\cal R}_{m+1}(g) &=& \zeta_m^{-1}\circ ({\cal R}_m(g))^{p_m} \circ \zeta_{m} \\
 & =& \left(\prod_{i=0}^m \zeta_i^{-1}\right)\circ g^{p_0\cdots p_m}\circ \left(\prod_{i=0}^k \zeta_i\right)
\end{eqnarray*}
where for every $m \geq 0$:
\begin{enumerate}
\item $\zeta_m$ is an accumulation point of the sequence $(\xi_{\phi(n)+m}(f))_n$,
\item ${\cal R}_{m+1}(g)$ is a limit point of the sequence  $({\cal R}_{\phi(n)+m+1}(f))_n$
\item $p_m$ is an accumulation point of the sequence $(a_{\phi(n)+m})_n\ . \ \  $
\cqd \end{enumerate}
\el

\begin{remark}
Let $g$ be an accumulation point of the renormalized maps $\RR_n(f)$. All along  this thesis, when we  say that $g$ is $(p_m)_{m \geq 0}$-infinitely renormalizable or we refer to the renormalized maps  $\RR_m(g)$, we always assume that $(\RR_m(g))_m$, their corresponding scaling maps $(\zeta_m)_m$ and $(p_m)_m$ satisfy statements 1, 2 and 3 of the above lemma.
\end{remark}

\noindent The renormalized maps $\RR_m(g)$ and their corresponding scaling maps $\zeta_m$ are respectively, limit points of the sequence of renormalized maps $(\RR_n(f))_n$ and the scaling maps $(\xi_n(f))_n$. Therefore,   the renormalized maps $\RR_m(g)$, the scaling maps  $\zeta_m$,  and their inverse $\zeta_m^{-1}$, are $C^{k+\epsilon}$  and their $C^i$-norm  for $1 \leq i \leq k$  are bounded by $M_f(i)$. Moreover,   $M_f(k)$ is also a  $\epsilon$-Holder constant  for $D^k\RR_m(g), D^k\zeta_m$ and $ D^k\zeta_m^{-1}(f)$.  This is summarized in the next result.
\bl \label{bondrenlema}
Any accumulation point $g$ of the sequence of renormalized maps $\RR_n(f)$ satisfies the bounded renormalization properties.
\el

\noindent The next lemma shows that, any accumulation point  of the sequence $\RR_n(f)$ can be seen  as a renormalized map of another accumulation point of the sequence
$({\cal R}_n(f))_n$.

\bl
For any accumulation point $g= \lim_{n\to \infty}{\cal R}_{\phi(n)}(f)$ of the sequence of renormalized maps of $f$, there is another limit point $ g_{-1}$ of the sequence $\RR_{\Phi(n)-1}(f)$  such that ${\cal R}(g_{-1}) \ = \ g$.
\el
\proof
We  choose a subsequence $\psi_1(n)$ of $\phi(n)$, such that
$({\cal R}_{\psi_1(n)-1}(f))_n$, $(\xi_{\psi_1(n)-1}(f))_n$ are convergent and $a_{\psi_1(n)-1}$ is constant. Let us define $g_{-1}=\lim_{n\to \infty}{\cal R}_{\psi_1(n)}(f)$, $\zeta_{-1} = \lim_{n\to \infty}\xi_{\psi_1(n)-1}(f)$ and $p_{-1}= a_{\psi_1(1)-1}$. Then we have
$${\cal R}(g_{-1})= \zeta_{-1}^{-1}\circ g_{-1}^{p_{-1}} \circ \zeta_{-1} =g \ . \ \ \Box$$

\noindent This lemma is the first step of an inductive process. For an arbitrary limit point  $g$ of the sequence $\RR_n(f)$ it gives another accumulation point $g_{-1}$ of the sequence $\RR_n(f)$ such that ${\cal R}(g_{-1}) \ = \ g$. For this $g_{-1}$ it gives an $g_{-2}$ such that ${\cal R}(g_{-2}) \ = \ g_{-1}$, and thus
${\cal R}_2(g_{-2}) \ = \ g$, In this way we obtain the following result:
\bl \label{g-l}
For any accumulation point $g$ of the sequence of renormalized maps $\RR_n(f)$,  there is a sequence of accumulation points $(g_{-\ell})_{\ell \geq 1}$ of the same sequence $\RR_n(f)$, such that
$${\cal R}_\ell(g_{-\ell}) \ = \ g \ . $$
\el

From now on, we will assume that $f$  satisfies the bounded
geometry hypothesis. This means  that  in addition to the bounded
renormalization properties, there exists constants $0<a_f < b_f
<1$  which depend only on
 ${f}$ such that for all $n\geq 0$, if ${\cal I}$ is an atom of the
 generation $l-1$ of $\RR_n(f)$ and $\JJ$, $\KK $ are  atoms of the generation $l$ of $\RR_n(f)$ with $\JJ, \KK \subset \II$, then
\begin{eqnarray}
a_f \ \leq \ {|{\cal J}| \over |{\cal I}|} \ \leq\ b_f & {\rm and }& a_f\ \leq \  {{\rm dist}({\cal J } ,{\cal K}) \over |{\cal I }|} \ \leq \ b_f \ . \label{IJKab}
\end{eqnarray}

\noindent The next result states that accumulation points of the sequence $\RR_n(f)$ also possess the bounded geometry property.

\begin{lemma} \label{boundgeoprol} Any accumulation point $g_\infty$  of the sequence $\RR_n(f)$ is an infinitely renormalizable map satisfying the bounded geometry properties. We can choose the bounds so that $a_{{g}} = a_f$ and $b_{{g}} = b_f$.
\end {lemma}
\proof From Lemma \ref{bondrenlema} we know that $g_\infty$ satisfies the bounded renormalization properties. It remains to prove that
given any atoms $\cal I$, $\cal J$, $\cal K$ of ${\cal R}_m(g_\infty)$,  with $\cal J$ and $\cal K$ of generation $l$, $\cal I$ of generation $l-1$ and $\cal J , K \subset I$. statement \ref{IJKab} holds.

\noindent However, since for every $m \geq 0$ the maps ${\cal R}_m(g_\infty)$  are also accumulation points of
$({\cal R}_{\phi(n)})_{n}$  we only need to prove that   any arbitrary accumulation point  $g$ of the sequence $\RR_n(f)$ satisfies statement \ref {IJKab} for any atoms $\cal I$, $\cal J$, $\cal K$ of $g$,  with $\cal J$ and $\cal K$ of generation $l$, $\cal I$ of generation $l-1$ and $\cal J , K \subset I$.

\noindent By definition,
\begin{eqnarray*}
{\cal I} &=&
g^{q({\cal I})} \circ \zeta_{0}(g) \circ  \cdots \circ\zeta_{l-1}(g)(\D^2) \\
{\cal J} &=&
g^{q({\cal J})} \circ \zeta_{0}(g) \circ  \cdots \circ\zeta_{l-2}(g)(\D^2) \\
{\cal K} &=&
g^{q({\cal K})} \circ \zeta_{0}(g) \circ  \cdots \circ\zeta_{l-2}(g)(\D^2)
\end{eqnarray*}
for integers $q({\cal I})$, $q({\cal J})$ and $q({\cal K})$ satisfying
\begin{eqnarray*}
0 \  \leq \ q({\cal I}) & \leq & p_0(g) \cdots p_{l-1}(g) -1 \\
0 \ \leq \ q({\cal J}) &\leq &p_0(g) \cdots p_{l-2}(g) -1 \\
0 \ \leq \ q({\cal K}) &\leq &p_0(g) \cdots p_{l-2}(g) -1  \ .
\end{eqnarray*}
Notice that the condition ${\cal J},{\cal K} \subset {\cal I}$  is  the same as  $q({\cal J})$ and $q({\cal K})$ being multiples of $q({\cal I})$.

\noindent Let us define,
\begin{eqnarray*}
{\cal I}_n&=& {\cal R}_{\psi(n)}^{q({\cal I})}(f) \circ \xi_{\phi(n)}(f) \circ \cdots \circ \xi_{\psi(n) +l-2}(f)(\D^2)\\
{\cal J}_n&=& {\cal R}_{\psi(n)}^{q({\cal J})}(f) \circ \xi_{\phi(n)}(f) \circ \cdots \circ \xi_{\psi(n) +l-2}(f)(\D^2)\\
{\cal K}_n&=& {\cal R}_{\psi(n)}^{q({\cal K})}(f) \circ \xi_{\phi(n)} (f)\circ \cdots \circ \xi_{\psi(n) +l-2}(f)(\D^2)
\end{eqnarray*}
where $\psi(n)$ is chosen so that
\begin{eqnarray*}
\lim_{n \to \infty}{\cal R}_{\psi(n)}(f) & = & g
\end{eqnarray*}
and
\begin{eqnarray*}
\lim_{n \to \infty}\xi_{\psi(n) + i}(f)&=& \zeta_i(g) \ \ \hbox{for $0 \leq i \leq l-1$}.
\end{eqnarray*}
For each $n$, ${\cal I}_n$, ${\cal J}_n$ and ${\cal K}_n$ are atoms of ${\cal R}_{\psi(n)}(f)$. We have that ${\cal J}_n$ and  ${\cal K}_n$ are atoms of generation $l$ and ${\cal I}_n$ is of generation $l-1$. Moreover, since $q({\cal I})$ divides $q({\cal J})$ and $q({\cal K})$ we have that ${\cal J}_n, {\cal K}_n \subset {\cal I}_n$.

\noindent From the bounded geometry of $f$ it follows that
\begin{eqnarray}
a_f \ \leq \ {|{\cal J}_n| \over |{\cal I}_n|} \ \leq\ b_f & {\rm and }& a_f\ \leq \  {{\rm dist}({\cal J }_n ,{\cal K}_n) \over |{\cal I }_n|} \ \leq \ b_f \ . \label{InJnKnab}
\end{eqnarray}

\noindent Since,  for any pair of subsets $A,B \subset \D^2$, the maps
$$
\begin{array}{ccl}
{\cal C}^0(\D^2)  & \longrightarrow &\R^+\\
h   & \longmapsto & |h(A)|
\end{array}
$$
and
$$
\begin{array}{ccl}
{\cal C}^0(\D^2)\times {\cal C}^0(\D^2) & \longrightarrow &\R^+\\
(h_1,h_2)  & \longmapsto & {\rm dist}(h_1(A),h_2(B))
\end{array}
$$
are continuous we get that
\begin{eqnarray*}
|{\cal I}| &=& \lim_{n \to \infty}|{\cal I}_n| \ ,\\
|{\cal J}| &=& \lim_{n \to \infty}|{\cal J}_n| \ ,\\
|{\cal K}| &=& \lim_{n \to \infty}|{\cal K}_n| \ ,\\
{\rm dist}({\cal J }, {\cal K }) &=& \lim_{n \to \infty}{\rm dist}({\cal J }_n,{\cal K }_n|) \ .
\end{eqnarray*}
Therefore these last four equalities  together with   statement \ref{InJnKnab} gives immediately that $\II$, $\JJ$, $\KK$ satisfy condition \ref{IJKab} and thus, the result is proved. \cqd

From now on we assume that, beyond the bounded geometry hypothesis,  the map $f$  contracts uniformly
the area, {\it i.e.} there exists $b$ such that
 $\vert \det(Df(x))\vert\leq b<1$ for all $x$
in $\D^2$. Then we have the following result:

\begin{lemma} \label{singular} Any  accumulation point  $g$ of the sequence $\RR_n(f)$ is a singular map,  i.e. $\det (Dg(x)) =  0$ for all $x\in \D^2$.
\end{lemma}
\proof We have
$$
{\cal R}_n(f) = \xi^{-1}_{n-1}(f)\circ \cdots \circ \xi^{-1}_0(f)\circ
f^{a_0 \cdots a_{n-1}} \circ \xi_0(f) \circ \dots  \circ
\xi_{n-1}(f).
$$
Since for any linear map $A$  in  dimension two, we have $|\det A|\leq
\Vert A\Vert^2$,  it follows that the modulus of the determinant of the scaling functions $\xi_i(f)$ and its inverse are bounded by $M_f^2(1)$, and thus:
$$ \vert \det D{\cal R}^n(f)(x)\vert \leq M_f^{4n}(1)\cdot b^{a_0 \cdots a_{n-1}} \leq M_f^{4n}(1)\cdot b^{2^n} \ ,
$$
and this quantity goes to 0 when $n$ goes to $\infty$. Thus, by continuity, any
accumulation point of the sequence $({\cal R}_n(f))_{n\geq 0}$ is a singular map. \cqd

\noindent Despite the fact that an accumulation point of the sequence of renormalized maps $\RR_n(f)$ is a singular map, the next lemma shows us that this map is not degenerate.

\begin{lemma} \label{xl} Let $g$ be a $C^1$ infinitely renormalizable map of the 2-disk. For each $l\geq 0$, there exists an atom $J_l$, of the
$l^{th}$ generation of $g$ and a point $x_l$ in $J_l$ such that $\Vert
Dg(x_l)\Vert \geq 1$.
\end{lemma}
\proof We know that
$$g^{p_0(g)\cdots p_{l-1}(g)}({\cal D}_l(g)) \subset {\cal D}_l(g),
$$
where
$$
{\cal D}_l(g) = \xi_{0}(g)\dots\xi_{l-1}(g)(\D^2),
$$
and that $g^{p_0(g)\cdots p_{l-1}(g)}$ permutes in ${\cal D}_l(g)$  the $p_l$ atoms of generation $l+1$ that belong to ${\cal D}_l(g)$. It follows that there exists a point $y_l$ in ${\cal
D}_l(g)$ such that  $$\Vert
Dg^{p_0(g)\cdots p_{l-1}(g)}(y_l)\Vert \geq 1.
$$
Consequently, in one of the ${p_0(g)\cdots p_{l-1}(g)}-1$ first images of
${\cal D}_l(g) $, that is to say, in  an atom $J_l$ of the
$l^{th}$ generation of $g$, there is a point $x_l$, image of $y_l$ by some
iterate of $f$,   such that $\Vert Dg(x_l)\Vert \geq 1$. \cqd

 Let $g$ be an arbitrary accumulation point  of  the
sequence of renormalized maps $\RR_n(f)$. By Lemma \ref{Rm(g)},
$g$ is $(p_m)_{m \geq 0}$-infinitely renormalizable and then Lemma
\ref{xl} gives a sequence  of points $(x_l)_{l \geq 0}$ such that
for every $l$, $x_l$ belongs to an atom $J_l$ of the $l^{th}$
generation  and $\Vert Dg(x_l)\Vert \geq 1$. We can extract from
it, a subsequence $(x_{\theta(l)})_l$ which converges to some $x
\in \D^2$.  By continuity, we get $\Vert Dg(x)\Vert \geq 1$. It
follows that there exists an open neighbourhood ${\cal V}_x$ of
$x$  such that, for all $y\in {\cal V}_x$, $\Vert Dg(y)\Vert \geq
1/2$. Thanks to Lemma \ref{singular}, we know that the map $g$ is
singular. Therefore, for all $y\in {\cal V}_x$, the dimension of
the kernel of $Dg$ is 1. Let us recall the following
trivialization theorem.

\bt  \label{foliationlemma}
let $U, V$ be  open sets of $\R^n$, $r$ an integer $\leq n$  and $F:U\to V$ a $C^k$ ($1 \leq k \leq \infty$) function such that $\dim(\ker DF(y) )=r$ for every $y\in U$.  Then for every $x\in U$ there exist two $C^k$ embeddings
$$u_{x,F}:[-1,1]^n\to U \  \ \ \hbox{and}   \  \  \   v_{x,F}:[-1,1]^n\to V$$
such that
$$
 u_{x,F}(0,\ldots ,0) \ =\ x \  \ \ \hbox{and}   \  \  \
  v_{x,F}(0,\ldots ,0) \ =\ F(x)
$$
and, denoting the points of $[-1,1]^n$ by $(t_1,t_2) \in [-1,1]^r \times  [-1,1]^{n-r}$  we have that
$$\begin{array}{cccll}
v_{x,F}^{-1}\circ F \circ u_{x,F} : & [-1,1]^r \times  [-1,1]^{n-r}& \longrightarrow & \R^n & \\
 & (t_1, t_2) & \longmapsto & (0,t_2) & \ .
\end{array}$$
\et

\begin{remark} \label{foliationremark}
 The maps $u=u_{x,F}, v=v_{x,F}$ in the above theorem depend
continuously on $F$ and $x$. In particular, the continuity with respect to $F$, ensures that if $(F_m)_m$ is a
$C^1$-converging  sequence of functions satisfying the hypothesis of
the theorem, then
$$\bigcap_{m \geq 0}{u_{x,F_m}([-1,1]^n)}$$
contains for some $\epsilon >0$ the ball $B_\epsilon(x)$.
\end{remark}

\begin{remark}
 For a map $F$ satisfying the hypothesis of Theorem \ref{foliationlemma}  we have a $C^k$
foliation given by
$u_{x,F}([-1,1]^r \times\{y\})$, for $y \in [-1,1]^{n-r}$. Since $v_{x,F}^{-1}\circ F \circ u_{x,F}
(t_1,t_2)=(0,t_2)$ it follows that $F\circ u(t_1,t_2) = v(0,t_2) = F\circ u(t_1,t_2) \ \forall -1 \leq t_1 \leq 1$. This means that $F$ maps the leaves
of the referred foliation into a point. These leaves are locally the leaves of the foliation given by
the field of directions defined by $\ker DF$.
\end{remark}

\noindent
 We can apply Theorem \ref{foliationlemma} to the map $g$. In this way, we obtain for the point $x \in V_x$, a change of coordinates $u_{x,g}:[-1,1]^2 \to V_x$ and  $v_{x,g}:[-1,1]^2 \to \R^2$ such that $u_{x,g}(0,0)=x$ and for every $y_1,y_2 \in [-1,1]$:
$$
g\circ u_{x,g}(y_1,y_2) \ = v_{x,g}(0,y_2) \ .
$$
This yelds that  for every integer $m$ and $y_1,y_2 \in [-1,1]$,
\begin{eqnarray}
g^m\circ u_{x,g}(y_1,y_2)&= & g^m\circ u_{x,g}(0,y_2) \ . \label{gmeq1}
\end{eqnarray}
From the fact that $u_{x,g}$ is an embedding mapping $(0,0)\in]-1,1[^2$ into $x$ we have that the point $x$ is in the interior of $u_{x,g}([-1,1]^2)$. Since $x=\lim_{l \to \infty}{x_{\theta(l)}}$, $x_{\theta(l)}$ will be in the interior of $u_{x,g}([-1,1]^2)$ for $l$ large.
Recall now that $g$ satisfies the bounded geometry properties  (Lemma \ref{boundgeoprol}) and then,  the diameter of
the atoms of the $l^{th}$ generation  goes to $0$ as $l$ goes to $\infty$.
Thus, there exists $l_0$ such that  for all $l \geq l_0, $
 $u_{x,g}([-1,1]^2)$  contains $J_{\theta(l)}$ (the atom of the $\theta(l)^{th}$
generation of $g$ containing $x_{\theta(l)}$). The map $g$ sends every  atom of the $\theta(l)^{th}$ generation  into
itself after $p_0(g)\cdots p_{\theta(l)}(g)$ iterations (we are using notation given by Lemma \ref{Rm(g)}):
\begin{eqnarray}
g^{p_0(g)\cdots p_{\theta(l)}(g)}(J_{\theta(l)}) &\subset&
J_{\theta(l)} . \label{g(J)}
\end{eqnarray}

\noindent By setting $I$ as the projection in the second coordinate of $u^{-1}_{x,g}(J_{\theta(l)})$, statement \ref{g(J)} allows us to consider the map
$g_{\theta(l)}\,=\,u_{x,g}^{-1}\circ\ g^{p_0(g)\cdots p_{\theta(l)}(g)}\circ u_{x,g}:
[-1,1]\times I \to\ [-1,1]\times I$. By \ref{gmeq1}  $g_{\theta(l)}$ reads:
$$
g_{\theta(l)}(x_1, x_2) \ \ = \ \ (g_{1,{\theta(l)}}(0,x_2), g_{2,
{\theta(l)}}(0,x_2)) \ .
$$

\noindent Notice that $I$ is a non-degenerate interval since $J_{\theta(l)}$ contains atoms of generation $p_0(g)\cdots p_{\theta(l)}(g)+1$ for $g$ (atom for $g^{p_0(g)\cdots p_{\theta(l)}(g)}$). Besides this, $I$ is closed, since it is a continuous image of the compact set $J_{\theta(l)}$. The change of coordinates $u_{x,g}:[-1,1]\times I \to W$ ($W=u_{x,g}([-1,1]\times I)$) gives rise, by  an affine map $\gamma:I\to [-1,1]$.  to a change of coordinates $\tau:[-1,1]^2\to W$ defined by
$$\tau(x_1,x_2) \ = \ u_{x,g}(x_1,\gamma^{-1}(x)) \ .$$
By definition, for $y \in [-1,1]$ the leaf $\tau([-1,1]\times\{y\}) \subset W$ intersects $J_{\theta(l)}$ in at least one point. Let $z$ be one point of the intersection  $J_{\theta(l)}\cap \tau([-1,1]\times\{y\})$. We have  from \ref{g(J)} that    $g^{p_0(g)\cdots p_{\theta(l)}(g)}(z) \in J_{\theta(l)} \subset W$. This means that
 $g^{p_0(g)\cdots p_{\theta(l)}(g)}(W)
\subset W$. Therefore we have:

\bl \label{reduct0} For any  accumulation point $g$ of the sequence of renormalized maps $\RR_n(f)$ we can find   an atom $J$ of the $l^{th}$ generation and a disk $W$ containing $J$ satisfying the following properties:
\begin{enumerate}
\item Denoting by  $q$ the returning time of $J$ into itself (by the map $g$) we have that $g^q(W) \subset W$.
\item  Under a $C^k$ change of coordinates $\tau: [-1,1]^2  \to W$, the map $H=\tau^{-1}\circ g^q\circ \tau$  reads
$$H(x_1,x_2) = (H_1(0,x_2), H_2(0,x_2)) \ .$$
\end{enumerate}
\el
\begin{figure}[!h]
\end{figure}

\noindent Setting $h(x)=H_2(0,x)$, $p_2(x_1,x_2) =x_2$ and $\iota_2(x)=(0,x)$, the above lemma implies that both diagrams
{\unitlength=1mm
\begin{center}
\begin{picture}(100,32)
\thinlines
\put(20,14){\makebox(0,0)[br]{$p_2\circ\tau^{-1}$}}
\multiput(24,20)(40,0){2}{\vector(0,-1){10}}
\put(68.5,14){\makebox(0,0)[bl]{$p_2\circ\tau^{-1}$}}
\multiput(22,22)(40,0){2}{\makebox(0,0)[bl]{W}}
\multiput(18,8)(40,0){2}{\makebox(0,0)[tl]{$[-1,1]$}}
\multiput(34,23)(0,-17.5){2}{\vector(1,0){20}}
\put(33.5,26){\makebox(0,0)[bl]{$g^{p_0(g)\cdots p_{\theta(l)}(g)}$}}
\put(43,2.5){\makebox(0,0)[tl]{$h$}}
\end{picture}
\end{center}
and
\begin{center}
\begin{picture}(100,32)
\thinlines
\put(20,14){\makebox(0,0)[br]{$\tau\circ\iota_2$}}
\multiput(24,20)(40,0){2}{\vector(0,-1){10}}
\put(68.5,14){\makebox(0,0)[bl]{$\tau\circ\iota_2$}}
\multiput(22,22)(40,0){2}{\makebox(0,0)[bl]{W}}
\multiput(18,8)(40,0){2}{\makebox(0,0)[tl]{$[-1,1]$}}
\multiput(34,23)(0,-17.5){2}{\vector(1,0){20}}
\put(43,26){\makebox(0,0)[bl]{$h$}}
\put(33.5,2.5){\makebox(0,0)[tl]{$g^{p_0(g)\cdots p_{\theta(l)}(g)}$}}
\end{picture}
\end{center}
} \noindent are commutative.  Thus, we can recover some of
the dynamics of $g^{p_0(g)\cdots p_{\theta(l)}(g)}$ in $W$,
by the dynamics of $h$.  In particular a nested sequence of
atoms ${(K_m)}_{m \geq 0}$ of $g$ of generation
$p_0(g)\cdots p_{\theta(l)}(g) + m$ in $W$ gives rise to a
cascade for $h$ defined by ${\cal D}_m(h) =
p_2\circ\tau(K_m)$ Therefore $h$ is infinitely
renormalizable and the atoms of generation $m$ of $h$ are
images by $p_2\circ\tau$ of atoms of generation
${\theta(l)}(g) + m$ of $g$ in $W$. The map $h$ will be
called a \em reduction map. \em Although the atom $K_m$ is
not necessarily an embedded disks (g is a degenerate map),
we can extend it to an embedded disk $\hat{K}_m$ by setting
$$\hat{K}_m \ = \ \tau^{-1}\left([-1,1]\times
p_2\circ\tau(K_m)\right)$$ With respect to this cascade we
have that $g^{p_0(g)\cdots p_{\theta(l)}(g)}|_W$ is
infinitely renormalizable.

We begin this section giving a precise meaning to the concept of reduction maps.

\bd \label{defreductmap} let $f$ be an infinitely renormalizable $C^k$ map of the two disk. We say that a $C^k$ map $h:[-1,1] \to [-1,1]$ is a \em reduction map \em for $f$ if we can find another map $h_1:[-1,1] \to [-1,1]$ and an $(p_m)_{m \geq 0}$-infinitely renormalizable map $g:\D^2\to \D^2$,   satisfying the following conditions:
\begin{enumerate}
\item There is a sequence of integers $\Phi(n)_{n \geq 0}$ such that
\begin{eqnarray*}
\lim_{n \to \infty}\RR_{\Phi(n)}(f) &\stackrel{C^k}{=}& g  .
\end{eqnarray*}
\item There exists  an embedded disk $W\subset \D^2$ containing an atom $J$  of the $l^{th}$ generation for $g$ such that $g^{q_l}(W) \subset W$, where  $q_l=p_0\cdots p_{l-1}$ is the returning time of $J$ to itself,
\item Under a  $C^k$ change of coordinates $\tau: [-1,1]^2 \to W$, $\tau^{-1}\circ g^{q_l} \circ \tau$ reads:
$$\tau^{-1}\circ g^{q_l}\circ \tau(x_1,x_2) = (h_1(x_2),h(x_2))  \ $$
\end{enumerate}
\ed

\begin{remark} \label{hredinfren}
Let ${(K_m)}_{m \geq 0}$ be a nested sequence of atoms in
$W$ for the map $g$,  with $K_m$ of generation $l+m$,
involved in the above definition. Then $h$ and $g^{q_l}|_W$
are $(p_{m+l})_{m \geq 0}$-infinitely renormalizable.  The
sequence $$\left(p_2\circ\tau(K_m)\right)_{m \geq 0}$$ is a
cascade for $h$ and $$\left(\tau^{-1}([-1,1]\times
p_2\circ\tau(K_m)\right)_{m \geq 0}$$ is a cascade for
$g^{q_l}|_W$.

\end{remark}

\begin{remark} \label{h-hprem}
For any atom $A$ of $h$,  the returning map $h^q_{|A}:A\to A$ is also a reduction map for $f$.
\end{remark}
\bt \label{reduct1a}
Let $1 \leq k \leq \infty$ and $f$ be an infinitely renormalizable $C^k$  map of the two disk which contracts the area and satisfies the bounded geometry hypothesis. Then $f$ admits a  reduction map.
\et
\proof Thanks to Corollary \ref{notempty}, we find an $C^k$-accumulation $g_0$ of the sequence of renormalized maps $(\RR_n(f))_n$ which is, by Lemma \ref{Rm(g)},  an infinitely renormalizable map.  Lemma \ref{reduct0} applied to $g_0$ gives immediately the result by setting  $h(x)=H_2(0,x)$ as the reduction map. \cqd

\noindent  By the previous theorem, reduction maps always exist  provided that $f$ contracts the area and satisfies the bounded geometry hypothesis.  We are going now to establish two results on the properties of  reduction maps for an infinitely renormalizable maps $f$.

\bt \label{reduct1} Let $f$ be an $(a_n)_{n \geq 0}$-infinitely renormalizable map satisfying the hypothesis of  Theorem \ref{reduct1a} and let $p_0$ be any given accumulation point of the sequence $(a_n)_{n \geq 0}$. Then   $f$ admits a reduction map which possesses a periodic point with period $p_0$.
\et
\proof: Let $g_0$ be the accumulation point of the sequence of renormalized maps  $\RR_n(f)$ given by Lemma \ref{notempty}. Choose a sequence  $\Phi_(n)$  of integers such that $a_{\Phi(n)}=p_0$, and $\lim_{n\to\infty}\RR_{\Phi(n)}(f)=g_0$. By Lemma \ref{gren313} we know that $g_0$ is renormalizable. More  precisely. $\RR(g_0)=\zeta_0^{-1}\circ g_0^{p_0}\circ\zeta_0$ for an embedding $\zeta_0$ which is an accumulation point of $\xi_{\Phi_(n)}(f)$. Since $g_0^{p_0}(\DD_0)\subset \DD_0$ for $\DD_0=\zeta_0(\D^2)$, it follows from the Brouwer's Fixed Point Theorem that  there is a point $y_0 \in \DD_0$ such that $g_0^{p_0}(y_0)=y_0$.
 Thanks to the bounded geometry hypothesis of $f$, we get from Lemma  \ref{boundgeoprol} that $g_0$ also satisfies the bounded geometry properties. In particular, since their distance is positive, the atoms of the first generation  are disjoint. This implies that $y_0$ is a periodic orbit with period $p_0$ for $g_0$.

\noindent Consider now the sequence $(g_{-l})_{l \geq 1}$ given by Lemma \ref{g-l}. Since $\RR_l(g_{-l})=g_0$, the existence of a periodic orbit $y_0$ for $g_0$ implies
that  in any atom of  generation $l$ for  the map $g_{-l}$, there is a periodic point  of period $p_0$ for the map $g_{-l}^{q_l}$, where $q_l$ denotes the returning time of any atom of generation $l$ of $g_{-l}$.

\noindent For every $l \geq 1$, Lemma \ref{xl}  gives  a point $x_{l}$  belonging to an atom $J_l$ of the $l^{th}$ generation  of $g_{-l}$ for which we have that $||Dg_{-l}(x_{-l})|| \geq 1$. From Theorem \ref{ascoli-arzela} the sequence $g_{-l}$  admits a converging subsequence in $\UU^k$. Denote by $(\theta(l))_{l \geq 1}$ a sequence of integers such that $g_{-\theta(l)}$ is convergent in $\UU^k$ and $x_{\theta(l)}$ also converges to a point $x \in \D^2$.
Then there exists $l_0$ such that, for $l \geq l_0$,  $||Dg_{-\theta(l)}(y)|| \geq 1/2$ for every $y$ in an open neighbourhood $V_x$ of $x$. This means, by Theorem \ref{singular} that $\dim(\ker(Dg_{-\theta(l)})) =1$ for every $y \in V_x$ and $l \geq l_0$.

\noindent  From now on, we will follow the same kind of  strategy as we have done to prove Lemma  \ref{reduct0}.

\noindent Thanks to Lemma \ref{foliationlemma} we obtain for $l \geq l_0$, an embedding $u_{x,g_{-\theta(l)}}:[-1,1]^2\to V_x$ such that by setting $g_{\theta(l)}\,=\,u_{x,g_{-\theta(l)}}^{-1}\circ\ g_{-\theta(l)}^{q_l}\circ u_{x,g_{-\theta(l)}}$ we have that
$$
g_{\theta(l)}(x_1, x_2) \ \ = \ \ (g_{1,{\theta(l)}}(0,x_2), g_{2,
{\theta(l)}}(0,x_2)) \ .
$$
Attending to Remark \ref{foliationremark} we have that for some $\epsilon >0$ and every $l \geq 1$, $u_{x,g_{-\theta(l)}}([-1,1]^2) \supset B_\epsilon(x)$. Applying Lemma \ref{boundgeoprol} to the map $g_{-\theta(l)}$ it follows that the atom $J_{\theta(l)}$ satisfies
\begin{eqnarray*}
|J_{\theta(l)}| &<& b_f^{\theta(l)-1}|\D^2| \ .
\end{eqnarray*}
Thus, for $l$ sufficiently large, $J_{\theta(l)} \subset B_\epsilon(x)$. Fix  an integer $l \geq l_0$ such that we have $J_{\theta(l)} \subset B_\epsilon(x)$.

\noindent Now, set  $W=u_{x,g}([-1,1]\times I)$, where $I$ is the projection in the second coordinate of $u^{-1}_{x,g_{-\theta(l)}}(J_{\theta(l)})$. By an affine map that carries $I$ to $[-1,1]$ we get from $u_{x,g_{-\theta(l)}}([-1,1]\times I)\to W$ a change of coordinates $\tau:[-1,1]^2\to W$.
It is straightforward to verify that $W \supset J_{\theta(l)}$ and $g_{-\theta(l)}^{q_{\theta(l)}}(W)\subset W$. If we set $g=g_{-\theta(l)}$, $J=J_{\theta(l)}$,  $H=\tau^{-1}\circ g_{-\theta(l)}^{q_{\theta(l)}} \circ \tau$ we have that
$$ \lim_{n\to \infty}{\cal R}_{\phi(n)}(f) \stackrel{C^k}{=} g \   $$
and
$$
H(x_1,x_2) = (H_1(0,x_2),H(0,x_2)) \ .
$$
Since the returning time of $J=J_{\theta(l)}$ into itself is $q_{\theta(l)}$, we have that the  map $h(x)=H_2(0,x)$ is a reduction map for $f$. Because of $g_{\theta(l)}^{q_{\theta(l)}}$ has a periodic orbit of period $p_0$ in $J_{\theta(l)}$, we have that $h$ has a periodic orbit of period $p_0$. \cqd

\bt \label{reduct2}
Let $f$ be an $C^\infty$ infinitely renormalizable  map of the two disk which contracts the area and satisfies the  bounded geometry hypothesis. Then  $f$ admits a reduction map $h$ that is multimodal.
\et
\proof From Lemma \ref{reduct2} we get that $f$ admits a $C^\infty$ reduction map $h$. From Remark \ref{hredinfren} we know that $h$ is infinitely renormalizable. Remark \ref{h-hprem} ensures that  any renormalized map of $h$  is also a reduction map  for $f$. Thus, to prove the result it suffices to show that there is some $n >0$ such that  $h$ has a finite number of critical points in each atom of the $n^{th}$ generation. We   prove this by contradiction.

\noindent Suppose that for all $n \geq 1$, $h$  has  in each generation $n$,  an atom $I_n$ with  infinitely many critical points. Denote by $z_n \in I_n$ an accumulation point of  critical points of $h$. In these conditions, $h'(z_n)=0$. Moreover, since $z_n$ is an accumulation point of infinitely many critical points and $h$ is $C^\infty$, $z_n$ is a flat point, i.e.  $h^{(k)}(z_n)=0$ for all $k \geq 1$.
From the bounded geometry of $g$ and attending to the fact that $p_2\circ\tau$
is a Lipschitz function ($\tau$ is the change of variables that appears in the Definition \ref{defreductmap})
any atom of generation $n$ has a diameter $d_n$ within the range
\begin{eqnarray}
{1 \over C }a^{p_0(g)\cdots p_{\theta(l)}(g) + n}_f \cdot |\D^2| \ \ \leq & d_n & \leq \ \ C\cdot b^{p_0(g)\cdots p_{\theta(l)}(g) + n}_f \cdot |\D^2| \nonumber \\
{1 \over C_1}\cdot a_f^n \ \ \leq & d_n & \leq \ \ C_1 \cdot b_f^n \ \ , \label{dn}
\end{eqnarray}
where $C_1$ is a constant that does not depend on $n$.

\noindent By   Taylor's expansion, for any $r \geq 1$ there exists $\delta(r) >0$ such that, for any interval $J$ in $I$ with diameter less than $\delta(r)$, we have
$$
|h(J)| \ \leq \ \sum_{k =1}^{r-1}{|h^{(k)}(z)||J|^k} +||h||_r |I|^r \ ,
$$
where $z$ is a point in $J$. So, for an integer $N=N(r)$ sufficiently large so that $C_1 \cdot b_f^n < \delta(r)$, for $n \geq N$, we obtain choosing  $z = z_n$, that the image of $I_n$ satisfies
\begin{eqnarray}
|h(I_n)| & \leq &   C_2(r) \cdot b_f^{rn} \ .\label{h(in)}
\end{eqnarray}
Since $h(I_n)$ is another atom of generation $n$, if we take an  $r$  such that $b_f^{r} < a_f$, we get for $n$ sufficiently large a contradiction between \ref{h(in)} and \ref{dn}. \cqd

\begin{remark}
From its proof we can see that it is possible to replace in  the previous Theorem the hypothesis on the  $C^\infty$ differentiability of $f$ by the weaker condition of  $f$ being $C^k$ for an $k$ satisfying $b_f^{k} < a_f$.
\end{remark}

\subsection{Proof of Theorem }
Since $f$ is a  $(2)_{n \geq 0}$-infinitely renormalizable
$C^\infty$ map that satisfies the bounded geometry
property, Theorem \ref{reduct2} says that $f$ has  a
reduction map $h:[-1,1]\to[-1,1]$ which is multimodal. The
fact that $h$ is a reduction map means that there is
another map $h_1:[-1,1] \to [-1,1]$ and an infinitely
renormalizable map $g:\D^2\to \D^2$,   satisfying the
following conditions:
\begin{enumerate}
\item There is a sequence of integers $\Phi(n)_{n \geq 0}$ such that
\begin{eqnarray*}
\lim_{n \to \infty}\RR_{\Phi(n)}(f) &\stackrel{C^\infty}{=}& g  .
\end{eqnarray*}
\item There exists  an embedded disk $W\subset \D^2$ containing an atom $J$  of the $l^{th}$ generation for $g$ such that $g^{q_l}(W) \subset W$, where  $q_l$  is the returning time of $J$ into itself,
\item Under a  $C^\infty$ change of coordinates $\tau: [-1,1]^2 \to W$, $\tau^{-1}\circ g^{q_l}(W) \circ \tau$ reads:
$$\tau^{-1}\circ g^{q_l}(W)\circ \tau(x_1,x_2) = (h_1(x_2),h(x_2))  \ \ .$$
\end{enumerate}
We notice that since $g$ is an accumulation point of the sequence of renormalized maps $\RR_n(f)$, we have by Theorem \ref{Rm(g)} that $g$ is $(2)_{m \geq0}$ infinitely renormalizable and $q_l=2^l$.

\noindent We extend the map $\tau$ to a homeomorphism $[-2,2]^2 \to \D^2$, and   continue to denote by $\tau$ this extension.

\noindent Let
 $(J_n)_{n \geq 0}$ be the sequence of disks, where $J_n$ is the  atom of generation $l$ of $\RR_{\phi(n)}$ such that $J_n$ converges in the Hausdorff topology to $J$.
The cascade of periodic orbits $(O_m(f))
_{m \geq 1}$ induces a cascade of periodic orbits $(O(\RR_k(f)))_{m \geq 1}$ for every renormalized map $\RR_k(f)$. By \ref{lnRk} the signature $s((O_m(f))_{m \geq 0})$ is an alternating sequence if and only if the signature $s((O(\RR_k))_{m \geq 1})$for any renormalized map  $\RR_k(f)$ alternates. Moreover, the map $\RR_{\phi(n)}^{2^l}:J_n\to J_n$ is conjugated to $\RR_{\phi(n)+2^l}$, and thus the corresponding signature is an alternating sequence if and only if the signature of $f$ is an alternating sequence.

\noindent To prove the Theorem \ref{mainobstruct} it suffices to find an integer $n$ such that the signature of the induced cascade of periodic orbits of $\RR_{\phi(n)}^{2^l}(f)$ in $J_n$ is an alternating sequence.

\noindent For every integer $n \geq1$, we define:
\begin{eqnarray*}
H_n    &=& \tau^{-1}\circ \RR_{\phi(n)}^{2^l}(f)\circ \tau\\
A_n    &=& \tau^{-1}(J_n) \ \ .
\end{eqnarray*}
\noindent  In $A_n$, the map $H_n$ has a cascade of periodic  orbits $(O_{n,m})_{m \geq 0}$ that is the conjugated cascade of periodic orbits of $\RR_{\phi(n)}^{2^l}$ in $J_n$. Thus, in order to prove Theorem \ref{mainobstruct} we are going to find an integer $n_0$ for  which the signature of $(O_{n_0,m})_{m \geq 0}$ is an alternating sequence.

\noindent We set also $H = \tau^{-1}\circ g^{2^l}\circ \tau$. For $m$ fixed, the sequence of periodic orbits of period $2^m$,  $(O_{n,m})_{n \geq
1}$ converges as $n\to\infty$ in the Hausdorff topology to an orbit $O_m$
of period $2^m$ for $H$,  contained in $[-1,1]^2$.
 The projection of $O_m$ by $p_2(x_1,x_2)=x_2$ is a periodic orbit
$O'_m$ for $h$. The sequence $(O'_m)_{m \geq 0}$ is a period doubling cascade  for $h$.

For each $H_n$, we can choose an
isotopy $(H_{n,t})_{t\in[0,1]}$, joining the identity $H_{n,0}$ to
$H_n=H_{n,1}$, in order  that the sequence of isotopies is
also convergent in the $C^0$-topology. It follows that   for each $m \geq 1$, there exists
$N(m)$ such that $l_m(H_n)$ is
constant for  $n \geq N(m)$. This allows us to extend the linking number $l_m$ to the degenerate map
$H$ by setting
\begin{eqnarray*}
 l_m(H) &=& \lim_{n \to \infty}l_m(H_n) \ .
\end{eqnarray*}
We call $({l_m(H) \over  2^m})_m$ the signature of the sequence
of periodic orbits $(O_m)_m$.

For an embedding $F$ of the disk, and two
distinct points $x,y \in {\D^2}$, denote by $\omega_F(x,y)$ the
algebraic number of half turns of the vector
$$(F_t(x)-F_t(y))/||F_t(x)-F_t(y)|| \ \ ,$$
where by half
turns we mean the (algebraic) crossings of the horizontal direction and
$(F_t)_{t\in[01]}$ is an isotopy between $F$ and the identity as before.

\noindent
Since the isotopies $(H_{n,t})_{t\in[0,1]}$ are convergent as $n \to \infty$, we have that for any sequence $(x_n,y_n)$ converging to $(x,y)$,  there
exists $N=N(x,y)$ such that $n \geq N$ implies that
$\omega_{H_n}(x_n,y_n)$ is constant. This allows us to define
\begin{eqnarray*}
\omega_H(x,y) &=& \lim_{n\to\infty}{\omega_{H_n}(x_n,y_n)} \ .
\end{eqnarray*}
  With this definition, we can write
\begin{eqnarray}
 l_m(H) &=& \sum_{i=0}^{2^m-1}{w_H(H^i(x),H^i(y))} \ , \label{l_n(H)=sum}
\end{eqnarray}
where $x$ is a periodic orbit of period $2^m$ which  is the
parent of the $2^{m+1}$-periodic orbit $y$.

\noindent
Now we recall that Theorem \ref{reduct2}   ensures that the reduction map $h$ is multimodal. Let $I_1,\ldots , I_q$ the intervals of monotonicity of $h$. The monotonicity
of $h$ in each $I_k$,  implies that $\omega_H$ is constant in $([-1,1]\times
I_k)^2$. We will denote by $M_k$ this constant.
\bl \label{lnphi}
We have that
\begin{eqnarray*}  l_m(H) &=& \sum_{k=0}^q 2M_k \Phi(k,m) \end{eqnarray*}
where $\Phi(k,m)$ denotes de number of points in $O_m$ which
belong to $[-1,1]\times I_k$, or equivalently, the number of points of $O'_m$ in $I_k$.
\el
\proof Pick  a point $x$ of $O_m$ in an atom of generation $n$ and a point $y$ in the same atom belonging to $O_{m+1}$.  Set $x_i=H^i(x)$ and $y_i=H^i(y)$, and denote  $p_2(t)$ by $t'$. We know that $x'_i$ is always between (its two sons) $y_i$ and $y_{i+2^m}$. This implies that
\begin{eqnarray*}
w_H(x_i,y_i)&=&w_H(x_i,y_{i+2^m})\ .
\end{eqnarray*}
From \ref{l_n(H)=sum} it follows that
\begin{eqnarray*}
 l_m(H) &=& \sum_{i=0}^{2^{m+1}-1}{w_H(x_i,y_i)} \\
    &=& \sum_{i=0}^{2^{m+1}-1}{w_H(x_i,z_i)} \ ,
\end{eqnarray*}
where $z_i=y_i$ if $x'_i,y'_i$ belong to the same $I_k$, or $z_i=y_{i+2^m}$ otherwise. This means that if $x'_i=x'_{i+2^m}
\in I_k$ then $w_H(x_i,z_i)=w_H(x_{i+2^m},z_{i+2^m})=M_k$. Therefore,
\begin{eqnarray*}
 l_m(H) &=& \sum_{i=0}^{2^{m+1}-1}{w_H(x_i,z_i)} \\
    &=& 2\sum_{i=0}^{2^m-1}{w_H(x_i,z_i)} \\
    &=& 2\sum_{k=0}^{q}M_k\Phi(k,m) \ . \ \ \Box
\end{eqnarray*}

\noindent
 From Proposition \ref{oliv} we know that
for all $m \geq 0$ and for all $k\in \{0, \ldots , q\}$
\begin{eqnarray*}  \Phi(k,m+1)&=&2\Phi(k,m)+r(k,m) \end{eqnarray*}
where for all $p_2 \geq p_1\geq 0$,
\begin{eqnarray*} \left|\sum_{j=p_1}^{p_2}{r(k,j)}\right| &\leq &2 \ . \end{eqnarray*}
 This together with  the previous lemma leads to:
\begin{eqnarray*}
 l_{m+1}(H)-2l_m(H) &=&
\sum_{k=1}^q{2M_kr(k,m)} \ \ .
\end{eqnarray*}
Therefore, by setting $M =\max_k{|M_k|}$ we have for arbitrary $p_2\geq p_1 \geq 0$,
\begin{eqnarray*}
\left|\sum_{m=p_1}^{p_2}l_{m+1}(H)-2l_m(H)\right| & \leq &
\sum_{k=1}^q\left|{2M_k}\left(\sum_{m=p_1}^{p_2}{r(k,m)}\right)\right|
\\
 & \leq & 4qM
\end{eqnarray*}
Thus,  since $l_{m+1}(H)-2l_m(H)$ are integers,  $l_{m+1}(H)-2l_m(H)$ has an
infinity of positive and an infinity of negative values.

\noindent
Then, for an arbitrary integer $N$ we choose $m_0$ and $m_1$ in order to have $l_{m_0+1}(H)-2l_{m_0}(H) >0$ and $l_{m_1+1}(H)-2l_{m_1}(H) <0$. But since $l_m(H)=\lim_{n\to \infty}{l_m(H_n)}$, there exists $n_0$, such that
\begin{eqnarray*}
  l_{m_0+1}(H_{n_0})& = &l_{m_0+1}(H) \\
 l_{m_0}(H_{n_0})& = &  l_{m_0}(H)\\
 l_{m_1+1}(H_{n_0})& = & l_{m_1+1}(H) \\
 l_{m_1+1}(H_{n_0})& = &  l_{m_1+1}(H) \ .
\end{eqnarray*}
Since $N$ is an arbitrary integer, we have that the signature of $(O_{n_0,m})_{m \geq 0}$ is an alternating sequence and consequently Theorem \ref{mainobstruct}    is proved.


\begin{thebibliography}{HH}

{\small



\bibitem {[Birm]} J. Birman, ``Braids, links and mapping class
group," {\sl Ann. of Math Stud. \bf 4} Princeton Univ. Press
(1984).

\bibitem {[BF]} R. Bowen and J. Franks, ``Invariant circles and the order structure of periodic orbits in monotone twist maps," {\sl Topology \bf 26} (1987) 21--35.

\bibitem {[BoHa]} P. Boyland and G. Hall, ``The periodic points of maps of the disk
and the interval," {\sl Topology \bf 15} (1976) 337--342.

\bibitem {[CaER]} M. Campanino, H Epstein and D. Ruelle, ``On the
existence of the Feigenbaum's fixed point," {\sl Comm. Math. Phys.
\bf 79} (1981) 261--302.

\bibitem {[CaE]} M. Campanino and H Epstein, ``On  Feigenbaum's
functional equation $g\circ g(\lambda x) + \lambda g(x)=0$," {\sl
Topology \bf 21} (1982) 125--129.

\bibitem {[Ca]} E. Catsigeras, ``Cascades of period doubling of stable codimension one," IMPA's thesis, 1994.

\bibitem {[CaGaMo]} E. Catsigeras, J-M Gambaudo and F. Moreira, ``Infinitely renormalizable diffeomorphisms of the disk at the boundary of chaos," to appear in {\sl
Proc. of Amer. Math. Soc. }

\bibitem {[CEK]} P. Collet, J-P Eckman and H. Koch, ``Period
doubling bifurcations for families of maps on $\R^n$" {J. Stat.
Physics, 25} (1980) 1-15.

\bibitem {[CT]} P. Coullet and C. Tresser, ``It\'erations d'endomorphismes et
groupe de renormalisation," {\sl J. Phys. \bf C5} (1978) 25--28.

\bibitem {[TC]}   P. Coullet and C. Tresser and ``It\'erations
d'endomorphismes et groupe de renormalisation," {\sl C. R. Acad. Sc.
Paris \bf287A} (1978) 577--580. }

\bibitem {Olivier} O. Courcelle, ``Cascades d'orbites
p\'eriodiques en dimension 1 et 2:hyperbolicit\'e et
renormalisaton", Thesis, Universit\'e de Nice-Sophia
Antipolis,1996.

\bibitem {[EW]} J-P. Eckman and P. Wittwer, ``A complete proof of
the Feigenbaum conjectures," {\sl J. Stat. Phys. \bf 46} (1987)
455--475.

\bibitem {[Fei]} M. Feigenbaum, ``A $C^2$ Kupka-Smale diffeomorphism of the
disk with no sorces or sinks," {\sl ``Dynamical Systems and Turbulence" -
Lecture Notes in Mathematics . \bf 898} (1980) 90--98.

\bibitem {[FY]} J. Franks and L.S, Young, ``Quantitative
universality for a class of non-linear transformations,"
{\sl J. Stat. Phys. \bf 19} (1978) 25--52.

\bibitem {[G]} J-M. Gambaudo, ``Boundary of Morse-Smale
Surface Diffeomorphisms: An Obstruction to Smoothness," in
{\ ``Proceedings of the International Conference on
Dynamical Systems and Related Topics"}, {\bf Advanced
Series in Dynamical systems, vol 9}  World Scientific,
London (1991) 141-152.


\bibitem {[GST]} J-M. Gambaudo, S. van Strien and C. Tresser,
``There exists a $C^\infty$ Kupka-Smale diffeomorphism on ${\bf S
}^2$ with neither sinks nor sources," {\sl Nonlinearity \bf 2}
(1989) 287--304.

\bibitem {[GST2]} J-M. Gambaudo, S. van Strien and C. Tresser,
``The periodic orbit structure of orientattion preserving
diffeomorphisms on $\D^2$ with topological entropy zero," {\sl
Ann. Inst. H. Poincar\'e \bf 49} (1989) 335--346.

\bibitem {[GSuT]} J-M. Gambaudo, D. Sullivan and C. Tresser,
``Infinite cascades of braids and smooth dynamical systems", {\sl
Topology \bf 33} (1994) 85--94.

\bibitem {[GT1]} J-M. Gambaudo and C. Tresser, ``How horseshoes are created,"
in {\bf `` Instabilities and Nonequilibrium Structures III"} E. Tirapegui and
W. Zeller Eds. (Reidel; Dordrecht/Boston/Lancaster/Tokyo) (1991).

\bibitem {[GT2]} J-M. Gambaudo and C. Tresser, ``Self-similar constructions in
smooth dynamics. Smoothness and dimension," {\sl Commun. Math. Phys. \bf 150}
(1992) 45--58.


\bibitem {[La1]} O.E. Lanford III, ``A computer -assisted proof of
the Feigenbaunm conjecure," {\sl Bull. Amer. Math. Soc.   \bf 6}
(1982) 427--434.

\bibitem {[MS]} W. de Melo and S. van Strien, {\bf One-Dimensional Dynamics}
(Ergebnisse der Mathematik und ihrer Grenzgebiete. 3. Folge, Vol.
25) (Spinger-Verlag, Berlin, 1993).



\bibitem {[Sm]} D. Smania, ``Complex bounds for multimodal maps: Bounded combinatorics,"
{\sl Nonlinearity. \bf 14 } (2001) 1311--1330.


\bibitem {[Su]} D. Sullivan, ``Bounds, Quadratic differentials, and
Renormalization Conjectures," in A.M.S. Centennial Publication, Vol. {\bf
(Providence, RI)} (1992).


\bigskip
\bigskip
\bigskip

\end{thebibliography}
\end{document}